\documentclass{article}
\usepackage[utf8]{inputenc}

\usepackage{graphicx}
\usepackage{amssymb}
\usepackage{soul}
\usepackage{amsmath}
\usepackage{amsfonts}
\usepackage{amstext}
\usepackage{amsbsy}
\usepackage{caption}
\usepackage{subcaption}
\usepackage{algpseudocode}
\usepackage{algorithm}

\usepackage{titlesec} 
\usepackage{abstract} 
\usepackage{titling} 
\graphicspath{{./}}

\usepackage{color}
\usepackage{xcolor}

\newcommand\fig[1]{Figure~\ref{#1}}
\newcommand\eq[1]{Eq.~\eqref{#1}}

\newcommand{\be}{\begin{equation}}
\newcommand{\ee}{\end{equation}}
\newcommand{\bs}[1]{\mathbf{#1}}
\DeclareMathOperator*{\argmin}{argmin}

\usepackage{lineno,hyperref}

\title{Parametric Dynamic Mode Decomposition for Reduced Order Modeling} 
\author{%
Quincy A. Huhn$^a$,
Mauricio E. Tano$^a$, 
Jean C.~Ragusa\thanks{Corresponding author}\ $^{,a}$, Youngsoo Choi$^b$ \\
\normalsize \href{mailto:quincy.huhn98@tamu.edu}{quincy.huhn98@tamu.edu}, 
\normalsize \href{mailto:mtano@tamu.edu}{mtano@tamu.edu}, 
\normalsize \href{mailto:jean.ragusa@tamu.edu}{jean.ragusa@tamu.edu}
\normalsize \href{mailto:choi15@llnl.gov}{choi15@llnl.gov}
\\
$^a$ Department of Nuclear Engineering\\Texas A\&M University, College Station, TX 77843\\
$^b$ Center for Applied Scientific Computing\\ Lawrence Livermore National Laboratory, Livermore, CA 94550
}
\date{} 
\renewcommand{\maketitlehookd}{%
\begin{abstract}
Dynamic Mode Decomposition (DMD) is a model-order reduction approach, whereby spatial modes of fixed temporal frequencies are extracted from numerical or experimental data sets. The DMD low-rank or reduced operator is typically obtained by singular value decomposition of the temporal data sets. For parameter-dependent models, as found in many multi-query applications such as uncertainty quantification or design optimization, the only parametric DMD technique developed was a stacked approach, with data sets at multiples parameter values were aggregated together, increasing the computational work needed to devise low-rank dynamical reduced-order models. In this paper, we present two novel approach to carry out parametric DMD: one based on the interpolation of the reduced-order DMD eigenpair and the other based on the interpolation of the reduced  DMD (Koopman) operator. Numerical results are presented for diffusion-dominated nonlinear dynamical problems, including a multiphysics radiative transfer example. All three parametric DMD approaches are compared.
\end{abstract}}

\begin{document}
\maketitle
\section{Introduction}
High-fidelity, dynamic, multiphysics simulations are becoming more routinely used in engineering practice and research.
Due to their larger computational cost, typically requiring several processing nodes or supercomputers, these simulations are most often limited to providing a single, best-estimate answer to a given problem~\cite{ashby2007}.
However, multi-query applications in engineering, such as uncertainty quantification, optimal control, or design optimization, require several, repeated evaluations of a computational model in which the uncertain or design parameters are assigned different values~\cite{bergmann2008, avramova2015, younis2018, huang2018, hu2018mp, zeng2019}.
Therefore, multi-query applications generally lie beyond the realm of high-fidelity simulation tools and rely on simplified, less-accurate computational models~\cite{sudret2017}.

Reduced order models (ROMs) offer the possibility of finding reduced representations of high-fidelity, dynamic, multiphysics models.
The ROM representations can significantly speed up computational times while introducing controllable prediction errors, when compared to the expected variations over the expected range of uncertain or design parameters in the multi-query problem~\cite{quarteroni2015, huang2017, tripathy2018}.
Broadly, ROMs can be classified as {\it intrusive} or {\it non-intrusive} regarding whether they need access or not to the discretized system of equations of the higher-fidelity simulations, respectively~\cite{chinesta2016}.
Both intrusive and non-intrusive ROMs involve the distillation of a reduced set of the the most observable and/or controllable modes from a set of snapshots obtained by parametric realizations of the high-fidelity system.
Typically, proper orthogonal decomposition (POD)~\cite{berkooz1993} or non-linear manifold learning via artificial neural networks (ANNs)~\cite{mohamed2018} are used for computing the reduced set of modes.
Then, intrusive ROMs use the computed modes for truncating, e.g., balanced truncation method~\cite{kurschner2018}, projecting, e.g., Galerkin or Petrov-Galerkin projections~\cite{benner2015, tr2021}, or moment-matching, e.g., moment-matching method~\cite{rafiq2021}, the high-fidelity discretized system of equations into a reduced representation.
Thanks to the temporal and parametric affinity of the high-fidelity model or by gappy-sampling of the nonlinear parametric terms, e.g., via the Discrete Empirical Interpolation Method~\cite{chaturantabut2010}, intrusive ROMs preserve temporal and parametric dependence and, thus, have been effectively used as surrogate models in multi-query problems~\cite{bazaz2012}.
Note, however, that in the dynamic evolution of the system, the most influential or observable modes will change in time and, thus, the projection performance may deteriorate during transients.
Nonetheless, variants for dynamically selecting the most influential modes have been introduced to overcome this issue in intrusive ROMs, such as spatio-temporal biorthogonal decomposition~\cite{aubry1991}, spectral POD~\cite{sieber2016}, time-dependent projection manifolds~\cite{loiseau2018}.
On the other hand, non-intrusive ROMs use the compressed set of modes and the dynamic evolution observed in the high-fidelity system realizations to derive a reduced expression for the evolution of the selected modes, see~\cite{xiao2015} as a review of the application of non-intrusive ROMs to the Navier Stokes equations.

As high-fidelity systems generally involve numerical artifacts for improving the computational performance or well-posedness of the system, e.g., matrix-free methods~\cite{adams2002} or artificial viscosity~\cite{guermond2011}, and are implemented on legacy code, the development and implementation of adapted, stable, intrusive ROMs require vast efforts~\cite{choi2021}.
Therefore, there is a growing interest in non-intrusive ROMs.
Dynamic mode decomposition (DMD) and ANNs have been widely used for developing non-intrusive ROMs.

Regarding ANNs, the reduction method generally relies on the reduction of high-fidelity snapshots into a latent reduced representation, via an encoder-type ANN or POD, coupled with an ANN expansion of the latent space into the high-fidelity space, via a decoder-type ANN~\cite{leecal2020}.
The parametric dependence of the high-fidelity system is emulated in ANN-driven ROMs via a parameter-dependent variational regularization of the latent reduced representation~\cite{mrosek2021, lee2021}.
Similarly, the time dependence is either emulated via variation regularization~\cite{eivazi2020} or is introduced in the ANN by adding recurrent structures between the encoder and decoder parts of the ANN~\cite{simpson2021}.
However, in their current state-of-the-art, ANN-driven ROMs are still subjected to numerical and performance constraints~\cite{fritzen2019, hesthaven2018, lui2019}.
Numerically, the modes selected by the ANN are non-orthonormal and involve representations with a lower order of convergence than the high fidelity model and that do not preserve physical constraints~\cite{otto2019}.
Regarding performance, the ANN compression and expansion involve a large number of weights between hidden layers, which frequently yield small speed-ups in terms of number of floating point operations.
Although improvements have been recently introduced to introduce physical constraints~\cite{erichson2019}, improve convergence order~\cite{pan2020}, and increase speed ups of ANN-driven ROMs~\cite{wu2020}, the development of efficient ANN-driven ROMs is still highly problem-dependent and requires a significant amount of user-tuning.

DMD methods have been popularly used for non-intrusive ROM over the past decade~\cite{kutz2016}.
DMD aims at finding a reduced representation Koopman operator~\cite{takeishi2017}, which allow us to naturally represent the dynamic evolution of the system on the selected POD modes via a transient, mode-decaying algebraic equation based on the modes and frequencies of this operator~\cite{arbabi2017}.
Although Mezic~\cite{mezic2005} was the first one to use DMD as a ROM method, many variants have later appeared to improve DMD performance on large-stream datasets~\cite{hemati2014}, prediction boundness~\cite{jovanovic2014}, ability to capture multiple scales~\cite{kutz2016b, dylewsky2019}, and reduce sensitivity to noise~\cite{dang2018}, between others.
The reader is referred to the review by Kutz et al. for details about these methods~\cite{kutz2016}.
Despite accounting for the system's time dependence, for a single parameter realization, comes naturally in DMD via the reduced representation of the Koopman operator, it is significantly more challenging to introduce parametric dependence on DMD-ROMs.

To the best of our knowledge, the only work proposing a parameter-dependent framework for DMD-ROM is the one by Syadi, Schmid et al.~\cite{sayadi2015}, where the DMD-ROM model shows a good performance for representing transient signals involving a bifurcation parameter. 
In the model proposed, the authors perform DMD on a virtual augmented snapshot matrix, created by time-stacking of high-fidelity snapshots for all parameter realizations, and then select the DMD modes by parametric-interpolation of the stacked DMD modes created.
Albeit performant in the problem evaluated, the proposed method presents three key limitations.
First, the computational cost of performing POD on the snapshot matrix grows at best linearly with the dimension of this matrix~\cite{singh2018}.
Thus, the memory trace and computational performance of this method may increase linearly with the number of parameters included in the training set of DMD, making it impractical for multi-query, multi-parametric applications.
Second, by assuming that snapshots can be vertically stacked, the method forces a parameter-independent frequency in the evolution of the DMD modes, which may be inadequate in nonlinear problems.
Finally, vertical stacking requires snapshots to be sampled at the same time, which introduces complications when time-adaptivity is used for efficiently solving the the high-fidelity models.
Hence, in this paper, we propose two new algorithms for performing parametric DMD-ROMs.
In these algorithms, DMD is performed independently per parameter realization and the parametric interpolation is performed at the level of the reduced Koopman operator.
Whereas one of this methods proposes the parametric-interpolation of the reduced Koopman operator eigenpair, the other one suggest an independent interpolation of each component of the Koopman operator.
We then compare the accuracy and computational performance of the method proposed by Syadi et al. against the two proposed methods, using 2D and 3D test cases. 

The rest of this article is organized as follows. In Section~\ref{sec:DMDbackground}, we discuss Dynamical Mode Decomposition (DMD) as a data-driven means to learn nonlinear dynamics, in the context of discretized governing laws arising from PDEs. In Section~\ref{sec:DMDparam}, we first review the state-of-the-art for parametric DMD (``stacked`` DMD) for parametric PDEs. We then propose two new and computationally more efficient approaches for parametric DMD, the reduced EigenPair Interpolation (rEPI) and the reduced Koopman Operator Interpolation (rKOI). In Section~\ref{sec:results}, we compare the various parametric DMD approaches for several parametric PDEs: (1) a nonlinear heat conduction problem, (2) an advection-reaction-diffusion problem, and (3) a coupled radiation energy/material temperature problem described by coupled nonlinear diffusion-reaction PDEs. We conclude in Section~\ref{sec:ccl}.

\section{Background on Dynamical Mode Decomposition}\label{sec:DMDbackground}
As a starting point, we consider a possibly nonlinear dynamical system:
\begin{equation}\label{eq:DynSyst}
\frac{\partial y}{\partial t} = \mathcal{F}(y(\bs{x},t),t) \,.
\end{equation}
Eq.~\eqref{eq:DynSyst} is obtained by invoking governing laws, typically expressed as partial differential equations. Examples of such systems are provided in Section~\ref{sec:results} for (1) nonlinear heat conduction, (2) advection-reaction-diffusion problem, (3) coupled temperature-radiation grey radiative transfer. 
In Eq.~\eqref{eq:DynSyst}, $t$ represents time, $\bs{x}$ is the phase-space independent variables (usually physical space coordinates, but other dimensions may be present, such as energy/frequency in non-grey radiative transfer, direction in particle transport, etc., when the solution belongs to a higher dimensional phase space), $u$ is the solution (state) of the governing law(s), and $\mathcal{F}(y,t)$ is a nonlinear operator that describes the governing law(s). If the physical problem contains several coupled physics, then $y$ is to be viewed as the solution of the whole multiphysics problem, where each component of $y$ is related to a single physics. Oftentimes, Eq.~\eqref{eq:DynSyst} depends on parameters that are either uncertain or that one may want to change or modify, so Eq.~\eqref{eq:DynSyst} may have to be solved repeatedly for each parameter realization. This is discussed further in Section~\ref{sec:DMDparam}.

After discretization of the phase space, one obtains a system of coupled ordinary differential equations (ODEs):
\begin{equation}\label{eq:DynSystDisc}
\frac{d \bs{y}}{d t} = \bs{F}(\bs{y}(t),t) \,,
\end{equation}
where $\bs{y} \in \mathbb{R}^N$ is a vector representing the state of the dynamical system for each of the $N$ degrees of freedom at time $t$ and $\bs{F}(\cdot)$ represents the dynamics in a discrete setting. The state vector $\bs{y}$ can be quite large, due to the large number of degrees of freedom involved in the discretization of the  phase space, that is, we often have $N \gg 1$. Finally, Eq.~\eqref{eq:DynSystDisc} is discretized in time and solution at various time steps $i\Delta t$ ($0 \le i \le m)$ are obtained.

In Dynamical Mode Decomposition (DMD), one attempts are replacing the operator $\bs{F}(\cdot)$ of Eq.~\eqref{eq:DynSystDisc} with a linear operator to capture the time evolution of the solution $\bs{y}$. Hence, we seek the following surrogate system
\begin{equation}\label{eq:DynSystDMD}
\frac{d\bs{y}}{d t} = \bs{\hat A} \bs{y} \,,
\end{equation}
which, in effect, assumes that the data collected was generated from linear dynamics and that linear operator $\bs{\hat A}$ approximates these dynamics. In DMD, the operator $\bs{\hat A}$ is learned in a data-driven fashion. We start with a time series of snapshots 
$\mathbf{S} = [\mathbf{y}_0, \mathbf{y}_1, \cdots, \mathbf{y}_m] = [\mathbf{y}_i]_{i=0}^m \in \mathbb{R}^{n \times m}$, which can be obtained from either experimental data or from simulations.
Here, $m$ stands for the number of snapshots in time and $n$ is the full-order dimension of the system, which can be either the number of degrees of freedom in numerical simulations, i.e., $n=N$, or the number of measuring devices in experiments.
We then split the data matrix $\mathbf{S}$ into a \textit{lagged matrix of snapshots} $\mathbf{S}^- = [\mathbf{y}_i]_{i=0}^{m-1} \in \mathbb{R}^{n \times (m-1)}$ and a \textit{forward matrix of snapshots} $\mathbf{S}^+ = [\mathbf{y}_i]_{i=1}^{m} \in \mathbb{R}^{n \times (m-1)}$.
In a continuous case, given the current state of a system, the operator that defines the new state of the system after an infinitesimal time is known as Koopman operator~\cite{mezic2015}.
Using this concept and the linearized version of the surrogate system in Eq.~\ref{eq:DynSystDMD}, one can define a discrete version of the Koopman operator $\mathbf{A} \in \mathbb{R}^{n \times n}$ as linear map from old to new states as follows:
\begin{equation}
    \label{Eq:KoopDef}
    \mathbf{S}^+ = \mathbf{A} \mathbf{S}^- \,.
\end{equation}
Note that $\bs{A}$ is the operator resulting from the time discretization of Eq.~\eqref{eq:DynSystDMD}. Obtaining $\bs{A}$ using data may seem a relatively simple, invoking the pseudo-inverse of the data as $\bs{A} = \bs{S}^+ (\bs{S}^-)^\dag$, where $^\dag$ denote the Moore-Penrose pseudo-inverse. However, the snapshot matrix can be large, usually with the number of (spatial) unknowns much larger than the number of collection times, i.e.,  $n \gg m-1$. Hence, $\bs{A}$ is a low-rank matrix, whose rank is at most $m-1$. Thus $\bs{A}$ is not sought after directly but rather we seek to discover its low-rank subspace first and use the dynamics of low-rank approximation to represent the full state dynamics, i.e., we seek for a reduced representation of the discrete Koopman operator $\bs{A}_r \in \mathbb{R}^{r \times r}$ with $ \ll n $. With DMD, one proposes to find this reduced representation based on the principal modes in space and time of the matrix of snapshots $\bs{S}$. To that effect, we compute the Singular Value Decomposition (SVD) of the lagged matrix of snapshots $\mathbf{S}^-$:
\begin{equation}
    \label{Eq:SVDS}
    \mathbf{S}^- = \mathbf{U} \mathbf{\Sigma} \mathbf{V}^T \,,
\end{equation}
where $\mathbf{U} \in \mathbb{R}^{n \times n}$ and $\mathbf{V}^T \in \mathbb{R}^{m \times m}$ are two orthogonal matrices that hierarchically contain the principal modes of the space and time behaviour of the system, respectively. The orthogonality of this matrix will be key to ensure that the algorithm will yield the correct modes for the reduced Koopman operator. Furthermore, $\mathbf{\Sigma} \in \mathbb{R}^{n \times m}$ is a matrix that stores the singular values, $\sigma_i$.

We then reduce the rank by truncating $\mathbf{U}$ and $\mathbf{V}$ to the first $r$ columns, i.e., defining $\mathbf{U}_r = [\mathbf{u}]_{i=1}^r \in \mathbb{R}^{n \times r}$ and $\mathbf{V}_r = [\mathbf{v}]_{i=1}^r \in \mathbb{R}^{m \times r}$, respectively. 
In addition, we also select the first $r$ singular values in $\mathbf{\Sigma}$ defining a diagonal matrix of singular values as $\mathbf{\Sigma}_r = \operatorname{diag}([\Sigma]_{i=1}^r) \in \mathbb{R}^{r \times r}$. The rank $r$ is typically chosen so that a certain fraction of the information (proportional to the singular values) is retained, that is,
\begin{equation} 
    \label{eq:rank_criterion}
    r = \argmin_{j} \frac{\sum_{i=1}^{i=j} \sigma_i}{\sum_{i=1}^{i=n} \sigma_i} < \tau \,,
\end{equation}
where $\tau$ is a user-specified fraction between 0 and 1 (usually close to 1). With the truncated SVD, we have
\begin{equation}
    \label{Eq:SVDRed}
    \mathbf{S}^- \approx \mathbf{U}_r \mathbf{\Sigma}_r \mathbf{V}_r^t \,.
\end{equation}
Then, using the orthogonality of $\mathbf{U}$ and $\mathbf{V}$ and plugging \eq{Eq:SVDRed} into \eq{Eq:KoopDef}, we obtain the reduced Koopman operator as
\begin{equation}\label{eq:koopmanop}
    \mathbf{A}_r \equiv \mathbf{U}_r^T \mathbf{A} \mathbf{U}_r = \mathbf{U}_r^T \mathbf{S}^+ \mathbf{V}_r \mathbf{\Sigma}_r^{-1} \,.
\end{equation}
$\mathbf{A}_r$ is of size $r \times r$, with $r$ significantly smaller than $n$ ($r \ll n)$, so $\mathbf{A}_r$ is a much smaller matrix than $\bf{A}$. It is, therefore, computationally inexpensive to perform the eigen-decomposition of $\mathbf{A}_r \in \mathbb{R}^{r \times r}$, yielding:
\begin{equation}\label{eq:redeigenmodes}
    \mathbf{A}_r \mathbf{W} = \mathbf{\Lambda} \mathbf{W} \,.
\end{equation}

The projected DMD modes are then computed as $\mathbf{\Phi} = \mathbf{U}_r \mathbf{W}$ and the evolution of the reconstructed system can then be computed via the formula
\begin{equation}
    \label{eq:reconstruction}
    \mathbf{y}(t) = \mathbf{\Phi} \mathbf{\Lambda}^{t/\Delta t} \mathbf{b}_0 = \sum_{i=1}^{r} b_{0i} \mathbf{\phi}_i (\lambda_i)^{t/\Delta t} \,,
\end{equation}
where $\bs{b}_0 = \mathbf{\Phi}^{\dagger} \mathbf{y}_0 \in \mathbb{R}^{n}$ is the vector of initial coefficients, $\mathbf{\Phi}^{\dagger}$ is the pseudo-inverse of the matrix of modes, and $\mathbf{y}_0 \in \mathbb{R}^{n}$ is the snapshot at the initial time. Algorithm \ref{alg:dmd} summarizes the steps for the classical DMD.

\begin{algorithm}
\caption{DMD Algorithm}\label{alg:dmd}
\begin{algorithmic}[1]
\State Solve Eq.~\eqref{eq:DynSystDisc} and collect temporal snapshots $\left[\mathbf{y}(t_i)\right]_{i=0}^m$
\State Arrange snapshots in $\mathbf{S}^+$ and $\mathbf{S}^-$ data matrices
\State Perform SVD of $\mathbf{S}^-$:   $\mathbf{S}^- = \bf{U \Sigma V}^T$
\State Retain $r$ modes and compute the reduced Koopman operator $\bf{A}_r$, \eq{eq:koopmanop}
\State Perform the eigen-decomposition of $\bf{A}_r$ to obtain the reduced eigenmodes, \eq{eq:redeigenmodes}
\State Recover the full-state modes $\bf{\Phi}=\mathbf{U}_r \mathbf{W}$
\State Reconstruct the full-state solution $\mathbf{y}(t)$, \eq{eq:reconstruction}
\end{algorithmic}
\end{algorithm}

\newpage

\section{Parametric Dynamical Mode Decomposition} \label{sec:DMDparam}
In this section, we expand upon the general algorithm for DMD, presented in Section~\ref{sec:DMDbackground}, in order to determine the temporal evolution of the system under parametric dependence. Indeed, the governing laws of Eq.~\eqref{eq:DynSyst} may be dependent upon some parameters $\mu$ and can be re-stated in a parametric fashion as follows:
\begin{equation}\label{eq:DynSystParam}
\frac{\partial y^\mu}{\partial t} = \mathcal{F}(y^\mu(\bs{x},t;\mu),t;\mu) \,.
\end{equation}
After discretization of the phase space, one obtains a system of coupled parametric ODEs:
\begin{equation}\label{eq:DynSystDiscParam}
\frac{d \bs{y}^\mu}{d t} = \bs{F}(\bs{y}^\mu(t;\mu),t;\mu) \,.
\end{equation}
In parametric DMD, one seeks to replace the operator $\bs{F}(\cdot)$ with a linear surrogate
\begin{equation}\label{eq:DynSystDMDParam}
\frac{d\bs{y}^\theta}{d t} = \bs{\hat A(\theta)} \bs{y}^\theta \,,
\end{equation}
where $\theta$ is a parameter realization (we used a different Greek letter to emphasize that the new parameter realization $\theta$ is not in the training set, which we denoted by using $\mu$'s).

To the best of our knowledge, the only parametric DMD algorithm is that of \cite{sayadi2015}. We will review that work as well as present two novel methods, developed in this paper. 
All three algorithms will be evaluated on time series from some parametric governing laws. The parametric DMD solutions will be assessed using  parameter values different than that of the training set.

\subsection{Stacked Parametric DMD}\label{sec:stackedDMD}
The current state-of-the-art in parametric DMD is the \textit{stacked} DMD algorithm, proposed \cite{sayadi2015}. This approach is very similar to the non-parametric DMD of Section~\ref{sec:DMDbackground}, except that the time series solutions for different parameter values are ``stacked'' to form an augmented snapshot matrix $\bs{S}$. The new snapshot matrix contains the vertically stacked time series for each parameter realization
\begin{equation}\label{eq:stackedSnaps}
    \mathbf{S} =\begin{bmatrix}
        \mathbf{S}_{\mu_1} \\
        \vdots\\
        \mathbf{S}_{\mu_i}\\
        \vdots\\
        \mathbf{S}_{\mu_{N_S}}
    \end{bmatrix} \,,
\end{equation}
where $N_S$ is the number of parametric realizations (i.e., samples) in the training set. Algorithm~\ref{alg:dmd} is then applied using this matrix of stacked snapshots ($\mathbf{S} \in \mathbb{R}^{(N_S\times n)\times m}$), yielding the parametric projected DMD modes $\mathbf{\Phi}$ as follows:
\begin{equation}\label{eq:stackedphi}
    \mathbf{\Phi} =\begin{bmatrix}
        \mathbf{\Phi}_{\mu_1} \\
        \vdots\\
        \mathbf{\Phi}_{\mu_i}\\
        \vdots\\
        \mathbf{\Phi}_{\mu_{N_S}}
    \end{bmatrix} \,.
\end{equation}
For any new parametric realization $\theta$, the projected DMD mode $\mathbf{\Phi}_\theta$ is obtained from a Lagrangian interpolation of a subset of the  $\mathbf{\Phi}_{\mu_j}$ vectors. That subset is defined as the set of training values that corresponds to the ball $\mathcal{B}$ of nearest neighbors for the new parameter $\theta$. That is, we pick parameters $\mu_j$'s from the training set such that $\mu_j \in \mathcal{B}(\theta)$. We define the set of such training parameters as: $\mathcal{S}=\left\{j\text{ such that } \mu_j \in \mathcal{B}(\theta)\right\}$ and thus use the corresponding components $\mathbf{\Phi}_{\mu_j}$ of $\mathbf{\Phi}$ ($j\in\mathcal{S}$) to obtain, by Lagrangian interpolation, $\mathbf{\Phi}_{\theta}$. A similar procedure is employed for the initial vector $\mathbf{b}_0^\theta$ that is interpolated from the $\mathbf{b}_0^{\mu_j}$ vectors.

In this approach, the DMD time eigenvalues are shared among all snapshots. Therefore, we can then simply use the time evaluation formula of Eq.~\eqref{eq:reconstructionParam} to produce the solution at the new parameter value $\theta$:
\begin{equation}
    \label{eq:reconstructionParam}
    \mathbf{y}^\theta(t) = \mathbf{\Phi}_\theta \mathbf{\Lambda}_\theta^{t/\Delta t} \mathbf{b}^\theta_0 = \sum_{i=1}^{r} b^\theta_{0i} \Phi_{\theta i} (\lambda_{\theta i})^{t/\Delta t} \,,
\end{equation}
Algorithm \ref{alg:dmdStacked} summarizes the steps needed for the stacked DMD.

\begin{algorithm}
\caption{Parametric DMD Algorithm via snapshot stacking}\label{alg:dmdStacked}
\begin{algorithmic}[1]
\State Solve Eq.~\eqref{eq:DynSystDiscParam} for all training parameters $\{\mu_j\}_{k=1}^{N_S}$ and collect temporal snapshots $\mathbf{S}_{\mu_j}=\left[\mathbf{y}(t_i; \mu_j)\right]_{i=0}^m$ for $1 \le j \le N_S$
\State Arrange the $\mathbf{S}_{\mu_j}$ matrices into the stacked matrix of snapshots $\mathbf{S}$, \eq{eq:stackedSnaps}
\State Arrange stacked snapshots in $\mathbf{S}^+$ and $\mathbf{S}^-$ data matrices
\State Perform SVD of $\mathbf{S}^-$:   $\mathbf{S}^- = \bf{U \Sigma V}^T$
\State Retain $r$ modes and compute the reduced Koopman operator $\bf{A}_r$, \eq{eq:koopmanop}
\State Perform the eigen-decomposition of $\bf{A}_r$ to obtain the reduced eigenmodes, \eq{eq:redeigenmodes}
\State Recover the full-state stacked DMD modes $\mathbf{\Phi}=\mathbf{U}_r \mathbf{W}$, \eq{eq:stackedphi}
\State Recover the initial coefficients for the DMD evolution as $\mathbf{b}_0 = \mathbf{\Phi}^\dag \left[\left(\mathbf{y}_0^{\mu_1}\right)^T \cdots \left(\mathbf{y}_0^{\mu_{N_S}}\right)^T \right]^T$
\State Interpolate $\mathbf{\Phi}$ and $\mathbf{b}_0$ to find the new parameter realization's modes ($\mathbf{\Phi}_\theta$) and initial condition ($\mathbf{b}^\theta_0$)
\State Reconstruct the full-state solution $\mathbf{y}^\theta(t)$, \eq{eq:reconstructionParam}
\end{algorithmic}
\end{algorithm}

\noindent
There are two limitations to the stacked DMD approach:
\begin{enumerate}
    \item The eigenvalues are shared among all parametric snapshot realizations, forcing the dynamics of the DMD reconstruction to be independent of the parametric realization.
    \item The computational cost and memory requirements of stacked DMD is large because SVD is performed on the data matrices containing all of the parametric snapshots. This can be prohibitively expensive for simulations with a large number of realizations of parameter values.
\end{enumerate}
 The two new methods we introduce below are meant to combat the limitations stacked DMD.

\subsection{Reduced Eigen-pair Interpolation (rEPI)}

The first method we propose is the reduced Eigen-Pair Interpolation (rEPI) method. Given a new parameter realization $\theta$, we only perform {\bf individual} classical Dynamical Mode Decompositions, one for each training parameter $\mu_j$ that is a nearest neighbor of the new parameter $\theta$ (i.e., $\mu_j \in \mathcal{B}(\theta$) or, equivalently $j \in \mathcal{S}$). We denote by $J=\text{card}(\mathcal{S})$ the number of nearest neighbors. Then, we perform a Lagrange interpolation of the eigenpairs of the reduced Koopman operator. Denoting the dimension of the input parameter space by $P$, the number of nearest neighbors is $J \propto 2^\text{dim}_\mu$ (e.g., as in the case of a tensor grid of the training set). So, the largest computational cost in this approach is the computation of $J$ SVD calculations, each one for a snapshot matrix of size $n \times (m-1)$ 
However, the important distinction with the stacked DMD approach as that we do not perform a unique SVD decomposition for the whole training set (stacked snapshot matrix of size $(P\times n) \times (m-1)$) but we perform individual SVD decompositions, one per training parameter. 
%
We then reduce the rank of these $J$ SVD decompositions according to the criterion in Eq~\eqref{eq:rank_criterion} and select the largest of the $J$ ranks ($r = \max_{j \in J}{(r_j)}$), yielding the reduced left- and right-singular vectors $\mathbf{U}_{jr}$ and $\mathbf{V}_{jr}^T$ as well as the reduced singular value matrix $\mathbf{\Sigma}_{jr}$. Now, can can construct the $J$ individual reduced Koopman operators
\begin{equation}\label{eq:koopmanop_indiv}
    \mathbf{A}_{jr} = \mathbf{U}_{jr}^T \mathbf{S}_j^+ \mathbf{V}_{jr} \mathbf{\Sigma}_{jr}^{-1} \,,
\end{equation}
from which we obtained the reduced eigenpairs:
\begin{equation}\label{eq:redeigenmodes_indiv}
    \mathbf{A}_{jr} \mathbf{W}_{jr} = \mathbf{\Lambda}_{jr} \mathbf{W}_{jr} \,.
\end{equation}
At this stage, the eigenpair of interest is obtained ($\mathbf{\Lambda}_{\theta r}$, $\mathbf{W}_{\theta r}$) via interpolation, using the $J$ eigenpairs ($\mathbf{\Lambda}_{jr}$, $\mathbf{W}_{jr}$).
Similarly, the interpolated SVD mode $\mathbf{U}_{\theta r}$ is computed via interpolation using the individual SVD modes $\mathbf{U}_{jr}$.
Then, the projected DMD modes are computed as $\mathbf{\Phi}_\theta = \mathbf{U}_{\theta r} \mathbf{W}_{\theta r}$. 
Finally, the initial value $\mathbf{b}_\theta$ is also interpolated using the $J$ neighbors $\mathbf{b}_j$ and the expression of Eq.~\eqref{eq:reconstructionParam} can be employed to reconstruct the time evolution of $\mathbf{y}^\theta$ as we now have $\mathbf{\Phi}_\theta$, $\mathbf{\Lambda}_{\theta r}$, and $\mathbf{b}_\theta$. The process is summarized in Algorithm~\ref{alg:dmdEigenPair}.

The rEPI method solves both limitations of stacked DMD for parametric applications in the most straightforward way. The time eigenvalues are no longer shared between all parameters as they are calculated as part of the DMD performed at each parameter value in the training set. Furthermore, since DMD is performed for each parameter independently of the others, the computational cost of this method scales linearly with the number of parameters used for the interpolation set ($J$) instead of the number of parameters used in the training set ($N_S$). Note, however, that when interpolating the eigenpairs, this method assumes smoothness of the found eigenpair along the parametric range. This limitation will limit the performance of this method for high-order parametric DMD as later observed in the results of this article.

\begin{algorithm}
\caption{Parametric DMD algorithm via reduced eigen-pair interpolation}\label{alg:dmdEigenPair}
\begin{algorithmic}[1]
\State Solve Eq.~\eqref{eq:DynSystDiscParam} for all training parameters $\{\mu_j\}_{k=1}^{N_S}$ and collect temporal snapshots $\bs{S}_{\mu_j}=\left[\mathbf{y}(t_i; \mu_j)\right]_{i=0}^m$ for $1 \le j \le N_S$
\State Arrange matrix of snapshots in $\mathbf{S}_j^+$ and $\mathbf{S}_j^-$ data matrices, $1\le j \le J$ 
\State Perform SVD of $\mathbf{S}_j^-$: $\mathbf{S}_j^- = \mathbf{U}_j \mathbf{\Sigma}_j \mathbf{V}_j^T$, $1\le j \le J$ 
\State Retain $r$ modes and compute the reduced Koopman operator $\mathbf{A}_{jr}$, , $1\le j \le J$, \eq{eq:koopmanop_indiv} 
\State Perform the eigen-decomposition of $\mathbf{A}_{jr}$ to obtain the reduced eigenmodes, \eq{eq:redeigenmodes_indiv}
\State Interpolate SVD-modes $\mathbf{U}_j$, eigen-modes $\mathbf{W}_{jr}$, and eigen-values $\mathbf{\Lambda}_{jr}$ for to find $\mathbf{U}_{\theta r}$, $\mathbf{W}_{\theta r}$, and $\mathbf{\Lambda}_{\theta r}$, respectively
\State Construct the DMD-mode $\mathbf{\Phi}_\theta = \mathbf{U}_{\theta r} \mathbf{W}_{\theta r}$
\State Recover the initial coefficients $\mathbf{b}_\theta$ by interpolating $\mathbf{b}_j$,  $1 \le j \le J$ 
\State Reconstruct the full-state solution $\mathbf{y}^\theta(t)$, \eq{eq:reconstructionParam}
\end{algorithmic}
\end{algorithm}

\subsection{Reduced Koopman Operator Interpolation (rKOI)}
In the second method we developed, we directly interpolate the reduced Koopman operator $\mathbf{A}_{\theta r}$ using the reduced operators $\mathbf{A}_{j r}$. That is, initially, the same procedure as in the rEPI method is followed. However, the reduced eigenpair  ($\mathbf{\Lambda}_{\theta r}$, $\mathbf{W}_{\theta}$)  is obtained directly from the eigenvalue problem 
$\mathbf{A}_{\theta r} \mathbf{W}_{\theta r} = \mathbf{\Lambda}_{\theta r} \mathbf{W}_{\theta r}$. The reduced Koopman operator $\mathbf{A}_{\theta r}$ is obtaining for a Lagrnagian interpolation of the $J$ nearest neighbors matrices, 
$\mathbf{A}_{j r}$.

The remainder of the procedure is then identical to the rEPI scheme: we interpolate $\mathbf{U}_{\theta r}$ in order to compute the projected DMD modes $\mathbf{\Phi}_\theta = \mathbf{U}_{\theta r} \mathbf{W}_{\theta r}$. We also interpolate the initial vector $\mathbf{b}_\theta$. Then, all the necessary information is at hand and Eq.~\eqref{eq:reconstructionParam} can be employed to reconstruct the time evolution of $\mathbf{y}^\theta$. 
The procedure is summarized as Algorithm~\ref{alg:dmdKO}.

The rKOI method also solves both limitations of stacked DMD. The time eigenvalues are no longer shared between all parameters as they are calculated as part of the eigen-decomposition of the interpolated reduced Koopman operator and this method scales linearly with the number of parameters in the interpolation set. Furthermore, note that this method assumes smoothness across the components of the reduced Koopman operator over the parametric range when performing interpolations. Since these components are solely a function of the projection of the matrix of snapthots over its SVD modes (Eq.~\eqref{eq:koopmanop}), this assumption is much less constraining that assuming smoothness over the reduced Koopman operator eigenpair.

\begin{algorithm}
\caption{Parametric DMD algorithm via reduced (discrete) Koopman operator interpolation}\label{alg:dmdKO}
\begin{algorithmic}[1]
\State Solve Eq.~\eqref{eq:DynSystDiscParam} for all training parameters $\{\mu_j\}_{k=1}^{N_S}$ and collect temporal snapshots $\bs{S}_{\mu_j}=\left[\mathbf{y}(t_i; \mu_j)\right]_{i=0}^m$ for $1 \le j \le N_S$
\State Arrange matrix of snapshots in $\mathbf{S}_j^+$ and $\mathbf{S}_j^-$ data matrices, $1 \le j \le J$ 
\State Perform SVD of $\mathbf{S}_j^-$: $\mathbf{S}_j^- = \mathbf{U}_j \mathbf{\Sigma}_j \mathbf{V}_j^T$, $1 \le j \le J$ 
\State Retain $r$ modes and compute the reduced Koopman operator $\mathbf{A}_{jr}$, $1 \le j \le J$, \eq{eq:koopmanop_indiv} 
\State Interpolate the components of $\mathbf{A}_{jr}$ for $\mu_j \in \mathcal{S}$ to find $\mathbf{A}_{\theta r}$
\State Perform the eigen-decomposition of $\mathbf{A}_{\theta r}$ to obtain the reduced eigen-pair $\left( \mathbf{\Lambda}_{\theta r} \mathbf{W}_{\theta r} \right)$
\State Interpolate SVD-modes $U_j$ to find $\mathbf{U}_{\theta r}$
\State Construct the DMD-mode $\mathbf{\Phi}_\theta = \mathbf{U}_{\theta r} \mathbf{W}_{\theta r}$
\State Recover the initial coefficients $\mathbf{b}_\theta$ by interpolating $\mathbf{b}_j$,  $1 \le j \le J$ 
\State Reconstruct the full-state solution $\mathbf{y}^\theta(t)$, \eq{eq:reconstructionParam}
\end{algorithmic}
\end{algorithm}

\section{Results}\label{sec:results} 
In this section, we test the accuracy and performance of the 3 parametric DMD methods presented above.
For test cases, we have selected problems that encompass various types of physical operators, including diffusion, advection, reaction, and coupled multi-physics problems. The problems are, nonetheless, diffusion-reaction dominated as it is well know that DMD for advection/wave propagation phenomena requires reorientation along the propagation characteristics \cite{LU2020109229}. 
The first two problems are variations on heat conduction, one being more diffusive and the other one more advection-dominated.
The last problem and the main focus of the results section is a coupled multi-physics radiative diffusion problem, where radiation energy and material temperature form a coupled system of differential equations that describes the absorption/re-emission of high-energy X-rays by matter.

The training set denotes the collection of snapshots, generated by sampling the parameter space (of dimension $P$). All tests are carried out for new realizations of the parameters, i.e., parameter values that are not in the training set. In order to assess the accuracy of the parametric DMD approaches, the DMD-reconstructed solution is compared against a reference solution $u^\text{ref}$, obtained by actually solving the full-order problem at the new parameter value.
To quantify the reconstruction errors from the various methods, we use the relative $L_2$ error norm, given in Eq.~\eqref{eq:rell2} where $u^\text{ref}$ is the reference snapshot from the testing set and $u^\text{DMD}$ is the DMD reconstructed solution at the parameter value. The relative $L_2$ error norm is calculated for each time step. We also provide the time-averaged relative $L_2$ error norm given in Eq.~\eqref{eq:timerell2-a} where $E$ is the relative $L_2$ error norm and $N_T$ is the number of time steps.
\begin{subequations}
\begin{equation}\label{eq:rell2}
    E(t_i,\mu_j) = \frac{\| u^\text{DMD}(t_i,\mu_j)-u^\text{ref}(t_i,\mu_j)\|_2}{\| u^\text{ref}(t_i)\|_2} \,,
\end{equation}
\begin{equation}\label{eq:timerell2-a}
     \bar E(\mu_j) = \frac{\sum_i^{N_T} E(t_i,\mu_j)}{N_T} \,.
\end{equation}
We also average these error norms over parameter values sampled from the testing set:
\begin{equation}\label{eq:timerell2-b}
     \langle E(t_i) \rangle = \frac{\sum_j^{N_S} E(t_i,\mu_j)}{N_S} \,,
\end{equation}
and
\begin{equation}\label{eq:timerell2-c}
      \langle \bar E  \rangle = \frac{\sum_j^{N_S} \tilde E(\mu_j)}{N_S} \,.
\end{equation}
\end{subequations}

\subsection{Transient Nonlinear Diffusion Problem}
\subsubsection{Problem Description}

\begin{figure}[!htpb]
    \begin{subfigure}{.33\linewidth}
         \centering
         \includegraphics[width=\textwidth]{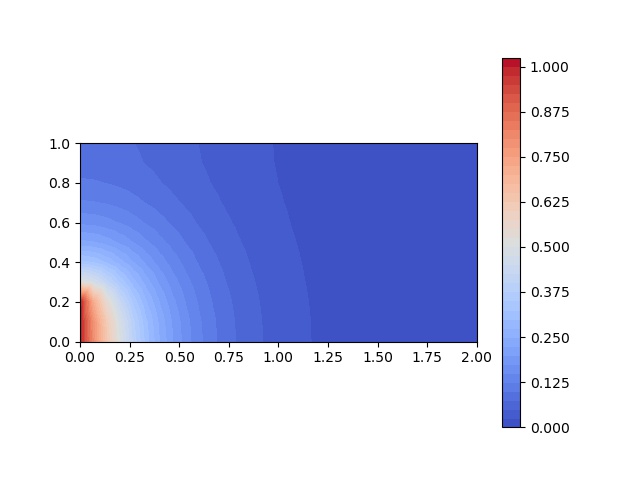}
         \caption{Solution at $b=0$, $t=0.1$}
         \label{fig:b1905}
    \end{subfigure}%
    \begin{subfigure}{.33\linewidth}
         \centering
         \includegraphics[width=\textwidth]{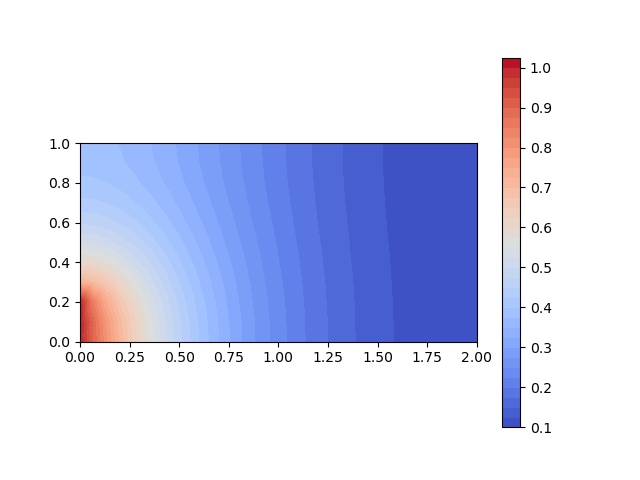}
         \caption{Solution at $b=0$, $t=0.5$}
         \label{fig:b191}
    \end{subfigure}
    \begin{subfigure}{.33\linewidth}
         \centering
         \includegraphics[width=\textwidth]{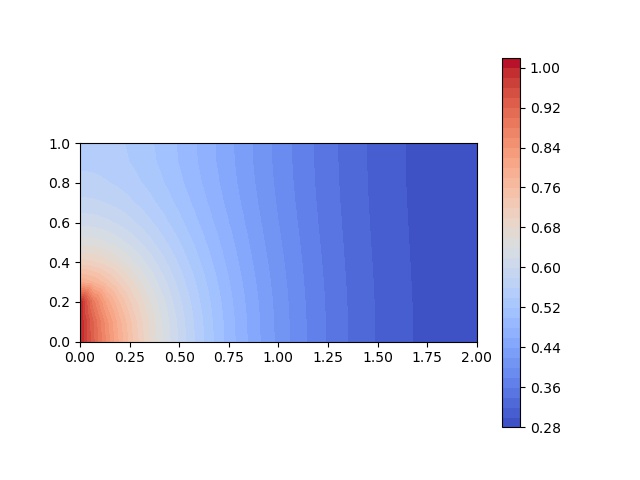}
         \caption{Solution at $b=0$, $t=1.0$}
         \label{fig:b192}
    \end{subfigure}
    \begin{subfigure}{.33\linewidth}
         \centering
         \includegraphics[width=\textwidth]{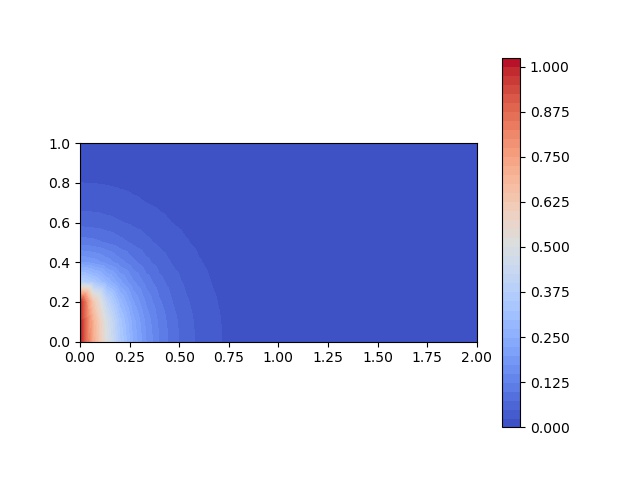}
         \caption{Solution at $b=4$, $t=0.1$}
         \label{fig:b2005}
    \end{subfigure}%
    \begin{subfigure}{.33\linewidth}
         \centering
         \includegraphics[width=\textwidth]{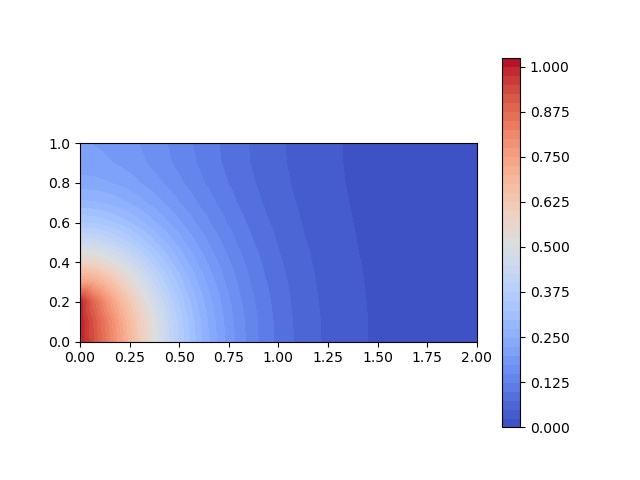}
         \caption{Solution at $b=4$, $t=0.5$}
         \label{fig:b201}
    \end{subfigure}
    \begin{subfigure}{.33\linewidth}
         \centering
         \includegraphics[width=\textwidth]{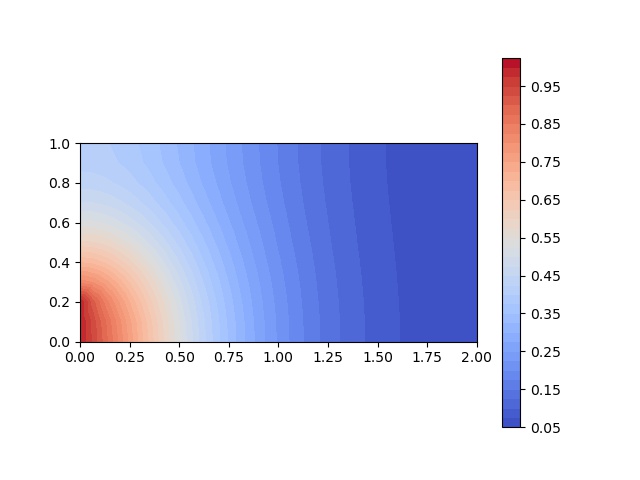}
         \caption{Solution at $b=4$, $t=1.0$}
         \label{fig:b202}
    \end{subfigure}
    \caption{Sample Snapshot Solutions of the Nonlinear Diffusion Problem at the Extremes of the Parameter $b$ Range Over Time.}
    \label{fig:nldiffsol}
\end{figure}

The first set of results deals with a transient nonlinear diffusion problem. The governing law is given in Eq.~\eqref{eq:nldiff}:
\begin{equation}\label{eq:nldiff}
    \frac{\partial T}{\partial t} + w \cdot \nabla T = \nabla \cdot k(T) \nabla T + f \qquad \forall \vec{r}\in \Omega
\end{equation}
where $\Omega = \left\{(x,y)\in [0,2] \otimes [0,1]\right\}$ is the problem domain, with a Dirichlet boundary condition $T(x=0,0>y<0.2,t)=1$ on the bottom of left wall, and homogeneous Neumann boundary condition elsewhere, $\nabla T \cdot n = 0$. 
The advection velocity for this problem is $w= 0.1\hat{i}+ 0\hat{j}$. The thermal conductivity $k$ is a non-linear function of temperature, given by Eq.~\eqref{eq:nonlineark}, where $a$ and $b$ are parameters of the conductivity correlation. The unknown parameter we vary in this problem is the exponent coefficient $b$, while $a$ is fixed to a value of $0.01$; the parameter range for $b$ is $0\leq b \leq 4$. Note that, when $b=0$, the problem is linear. Sample snapshot solutions for $b=0$ and $b=4$  are given in Fig.~\ref{fig:nldiffsol}, at different time instants.
\begin{equation}\label{eq:nonlineark}
    k(T) = a + T^b \,.
\end{equation}

\clearpage
\subsubsection{Results}
\begin{figure*}
    \centering
    \begin{subfigure}[b]{0.475\textwidth}
        \centering
        \includegraphics[width=\textwidth]{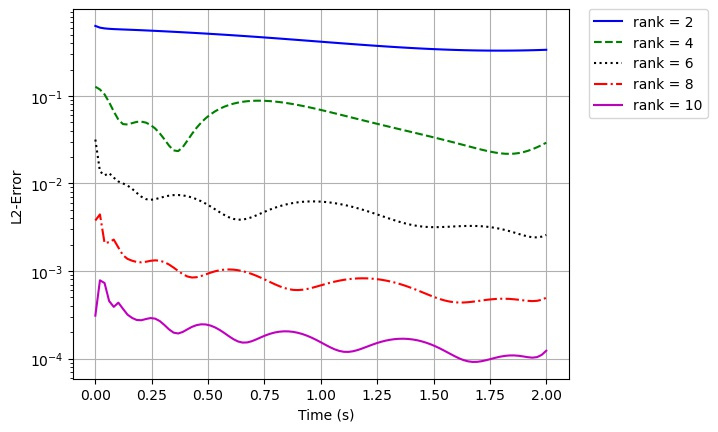}
         \caption{Stacked DMD}
         \label{fig:stackedfront}
    \end{subfigure}
    \hfill
    \begin{subfigure}[b]{0.475\textwidth}  
        \centering 
        \includegraphics[width=\textwidth]{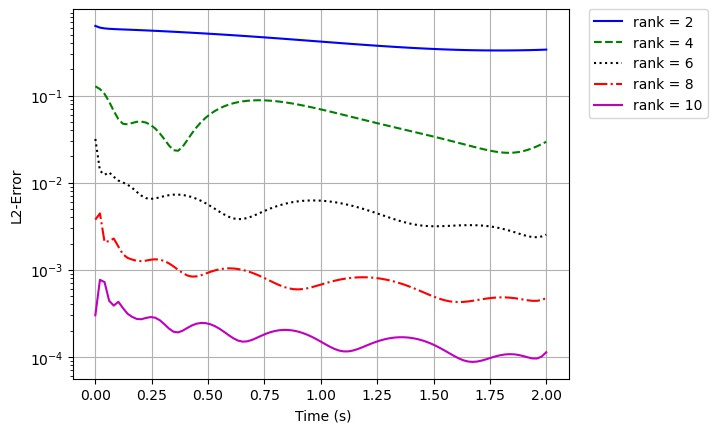}
        \caption{rKOI}
        \label{fig:kofront}
    \end{subfigure}
    \vskip\baselineskip
    \begin{subfigure}[b]{0.475\textwidth}   
        \centering 
        \includegraphics[width=\textwidth]{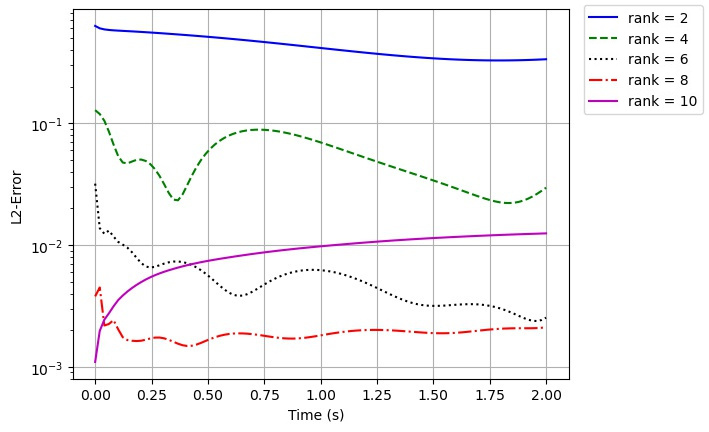}
        \caption{rEPI}
        \label{fig:epairmfront}
    \end{subfigure}
    \hfill
    \caption{Parameter-averaged, Time-dependent Relative $L_2$ Error, $\langle E(t_i) \rangle$, as a function of Rank for the Three Parametric DMD Methods on the Nonlinear Diffusion Problem.}
    \label{fig:compareFront}
\end{figure*}

In this problem, we generated 50 snapshots with varying thermal conductivity non-linearity exponents $b$ equally spaced in $[0,4]$. A testing set of 20 snapshots, sampled within the parameter range, was used. We compute the relative $L_2$ error norm over the testing set, $\langle E(t_i) \rangle$. Data points from the testing set are added until the change in the $\langle \bar E  \rangle$ error is less than $0.1\%$ or until all of the testing set samples have been utilized. This process is employed in all subsequent results as well, unless stated otherwise. In Fig.~\ref{fig:compareFront} we show the $\langle E(t_i) \rangle$ error and in Table~\ref{tab:aderror} we show the $\langle \bar E  \rangle$ error of each of the three methods as a function of DMD reconstruction rank.

\begin{table}[]
    \centering
        \begin{tabular}{|l|l|l|l|l|l|}
        \hline
        Rank                 & 2   & 4   & 6   & 8   & 10\\
        \hline
        Stacked Error & 4.31 E-1 & 5.21 E-2 & 5.48 E-3 & 8.86 E-4 & 1.89 E-4 \\
        \hline
        rKOI Error    & 4.37 E-1 & 5.21 E-2 & 5.57 E-3 & 8.82 E-4 & 1.87 E-4\\
        \hline
        rEPI Error    & 4.34 E-1 & 5.24 E-2 & 5.51 E-3 & 7.91 E-3 & 2.44 E-2 \\
        \hline
        \end{tabular}
        \caption{$\langle \bar E  \rangle$ Error in Rank for the Three Parametric DMD Methods on the Nonlinear Diffusion Problem}
        \label{tab:aderror}
\end{table}

For this problem, the error monotonically decreases as a function of the approximation rank $r$. This pattern will hold on all problems for most of the methods. On this problem, stacked DMD and reduced Koopman operator interpolation perform very similarly while the reduced Eigen-pair interpolation method presents higher errors. The three methods perform the same for rank 2, 4, and 6 reconstructions but rEPI has slightly worse error for rank 8 reconstruction and rank 10 reconstruction breaks with the general trend by having an increase in error with a higher rank reconstruction. This instability in rank is an undesirable trait because it is impossible to know a priori if increasing rank will decrease the error, making the rEPI method of questionable use for parametric DMD.

\subsection{Incident-Jet Problem}

\subsubsection{Problem Description}
\begin{figure}
    \begin{subfigure}{.33\linewidth}
         \centering
         \includegraphics[width=\textwidth]{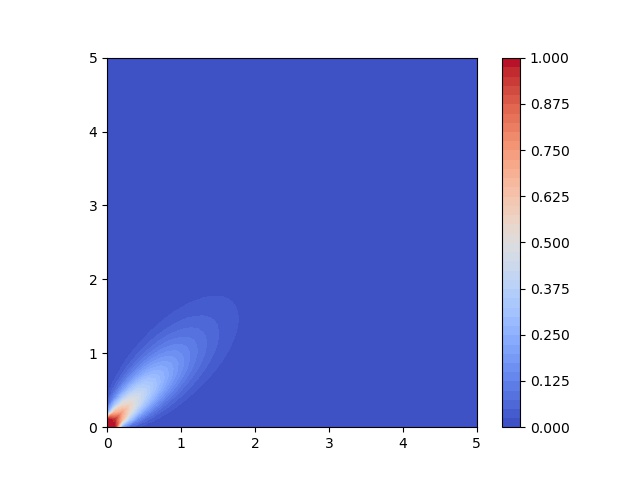}
         \caption{$k=0.2$, $t=0.2$}
         \label{fig:k021}
    \end{subfigure}%
    \begin{subfigure}{.33\linewidth}
         \centering
         \includegraphics[width=\textwidth]{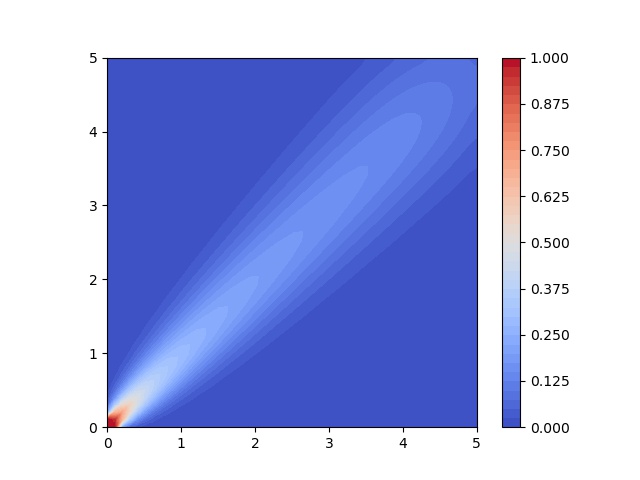}
         \caption{$k=0.2$, $t=1.0$}
         \label{fig:k022}
    \end{subfigure}
    \begin{subfigure}{.33\linewidth}
         \centering
         \includegraphics[width=\textwidth]{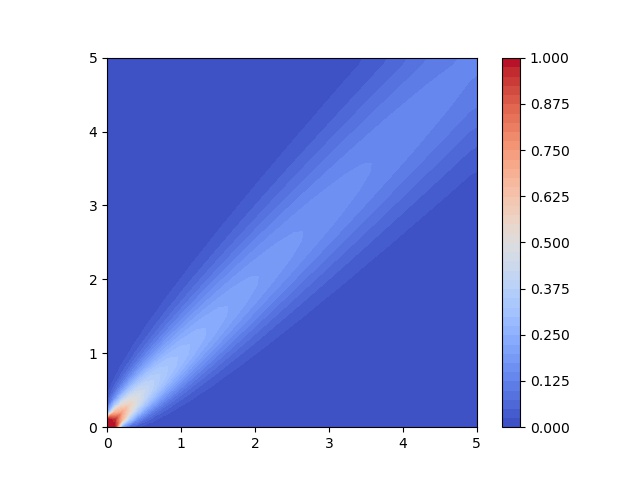}
         \caption{$k=0.2$, $t=2.0$}
         \label{fig:k024}
    \end{subfigure}
    \begin{subfigure}{.33\linewidth}
         \centering
         \includegraphics[width=\textwidth]{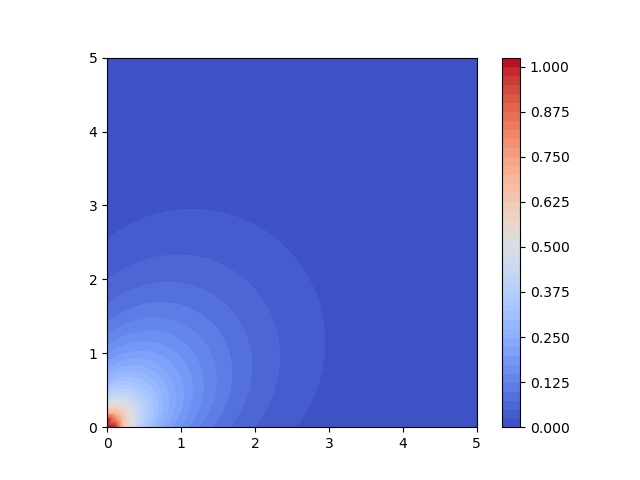}
         \caption{$k=5.0$, $t=0.2$}
         \label{fig:k031}
    \end{subfigure}%
    \begin{subfigure}{.33\linewidth}
         \centering
         \includegraphics[width=\textwidth]{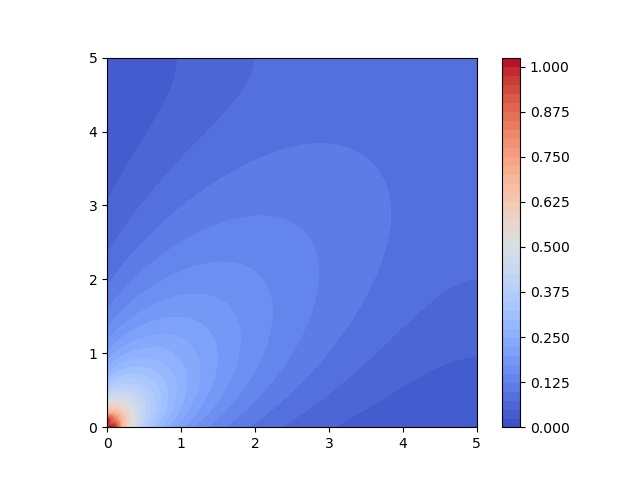}
         \caption{$k=5.0$, $t=1.0$}
         \label{fig:k032}
    \end{subfigure}
    \begin{subfigure}{.33\linewidth}
         \centering
         \includegraphics[width=\textwidth]{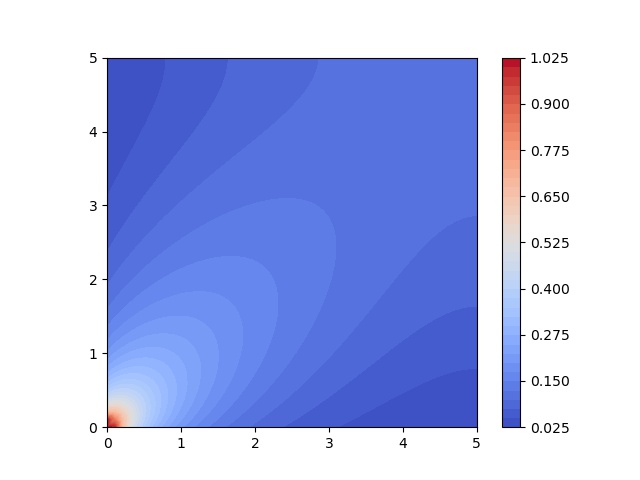}
         \caption{$k=5.0$, $t=2.0$}
         \label{fig:k034}
    \end{subfigure}
    \caption{Sample Solutions of the Incident Jet Problem at the Extremes of the Parameter Range Over Time}
    \label{fig:jetsol}
\end{figure}

The second problem is an incident-jet problem. Compared to the previous nonlinear diffusion problem, the incident-jet test case is more advection-dominated. The thermal conductivity was chosen to be independent of temperature and varies in $0.2\leq k \leq 5$. The governing law is given in Eq.~\eqref{eq:jet}. 
\begin{equation}\label{eq:jet}
    \frac{\partial T}{\partial t} + w \cdot \nabla T = \nabla \cdot k \nabla T
\end{equation}
The domain for this problem is $\Omega = \left\{(x,y)\in[0,5]\otimes [0,5]\right\}$ with a Dirichlet boundary condition of $T(x<0.1,y<0.1,t)=sin(10\pi t)$ in the bottom left corner, and a homogeneous Neumann boundary condition of $\nabla T \cdot n = 0$ elsewhere. In this problem, the advection velocity $w$ has $x$ and $y$ components such that $w = 5\hat{i} + 5\hat{j}$. The parameter we chose to vary in this problem is the thermal conductivity, the parameter range for $k$ is $0.2\leq k \leq 5$. Example solutions at $k=0.2$ and $k=5$ over time are given in Fig. \ref{fig:jetsol}. We compare the performance of the three parametric DMD methods on this problem in the following section.

\subsubsection{Problem Results}
\begin{figure*}
    \centering
    \begin{subfigure}[b]{0.475\textwidth}
        \centering
         \includegraphics[width=\textwidth]{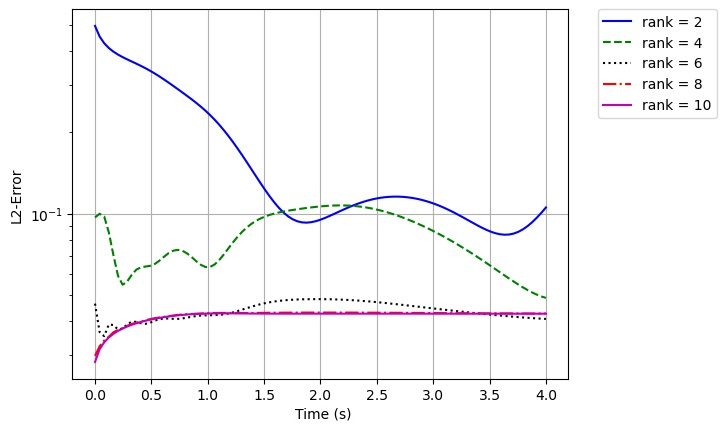}
         \caption{Stacked DMD}
         \label{fig:stackedjet}
    \end{subfigure}
    \hfill
    \begin{subfigure}[b]{0.475\textwidth}  
        \centering 
        \includegraphics[width=\textwidth]{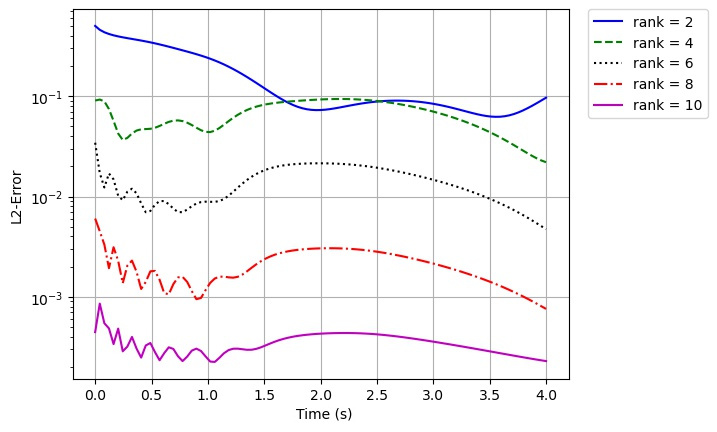}
         \caption{rKOI}
         \label{fig:kojet}
    \end{subfigure}
    \vskip\baselineskip
    \begin{subfigure}[b]{0.475\textwidth}   
        \centering 
        \includegraphics[width=\textwidth]{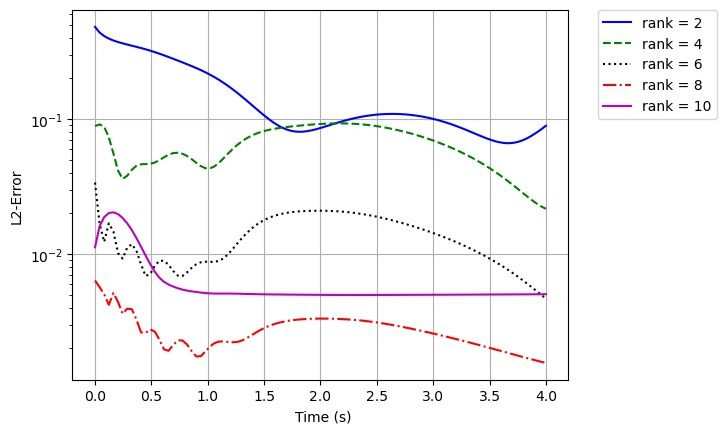}
         \caption{rEPI}
         \label{fig:epairmjet}
    \end{subfigure}
    \hfill
    \caption{$\langle E(t_i) \rangle$ Error as a Function of Rank for the Three Parametric DMD Methods on the Incident-Jet Problem}
    \label{fig:comparejet}
\end{figure*}

In this problem, we collected snapshots for 50 thermal conductivities $k$ equally spaced in $[0.2,5]$. A testing set was created by randomly setting aside 20 snapshots, while the training set contained the remainder snapshots. To get a relative $L_2$ error norm over the testing set we calculated both the $\langle E(t_i) \rangle$ error and the $\langle \bar E  \rangle$ error.  In Fig. \ref{fig:comparejet} we show the ensemble-averaged error $\langle E(t_i) \rangle$ and in Table~\ref{tab:jeterror} we show the ensemble- and time-averaged error $\langle \bar E  \rangle$ of each of the three methods, as a function of DMD reconstruction rank.

\begin{table}[]
    \centering
        \begin{tabular}{|l|l|l|l|l|l|}
        \hline
        Rank                 & 2   & 4   & 6   & 8   & 10\\
        \hline
        Stacked Error & 1.67 E-1 & 8.45 E-2 & 5.21 E-2 & 5.00 E-2 & 4.99 E-2 \\
        \hline
        rKOI Error    & 1.59 E-1 & 6.57 E-2 & 1.41 E-2 & 2.14 E-3 & 3.66 E-4\\
        \hline
        rEPI Error    & 1.60 E-1 & 6.52 E-2 & 1.40 E-2 & 2.47 E-3 & 3.64 E-3 \\
        \hline
        \end{tabular}
        \caption{$\langle \bar E  \rangle$ as a Function of Rank for the Three Parametric DMD Methods on the Incident-Jet Problem}
        \label{tab:jeterror}
\end{table}

In this problem rKOI is the clear winner out of the three methods. The stacked DMD error saturates at a relatively large value, even when rank is increased. This may be due to the fact that DMD time eigenvalues are generated and shared among all training snapshots, a drawback that is emphasized more in this advection-dominated problem. Because this method is the state of the art, we will continue employ it as baseline comparison in the remainder of this paper, even though its performance may be low at times. The rEPI method exhibits the same behavior as in the previous transient nonlinear diffusion problem, where the error does not monotonically decrease with rank. This behavior is disqualifying for a parametric DMD method due to the unpredictability of its behavior in rank. The rKOI method has also performed better than the rEPI approach for two test problems so far, thus we conclude that there are little benefits in the rEPI approach. As such, we will no longer present results using the rEPI approach.

\subsection{Radiative Diffusion Problem}

\subsubsection{Problem Description}
\begin{figure}
    \begin{subfigure}{.5\linewidth}
         \centering
         \includegraphics[width=\textwidth]{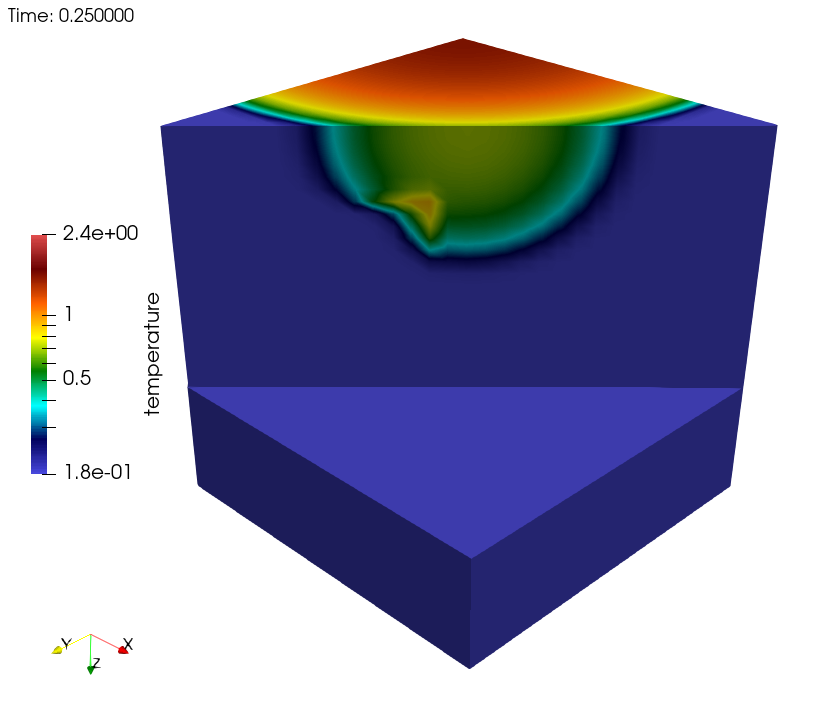}
         \caption{Solution at $Z=10$, $t=0.25$}
         \label{fig:z10t025}
    \end{subfigure}%
    \begin{subfigure}{.5\linewidth}
         \centering
         \includegraphics[width=\textwidth]{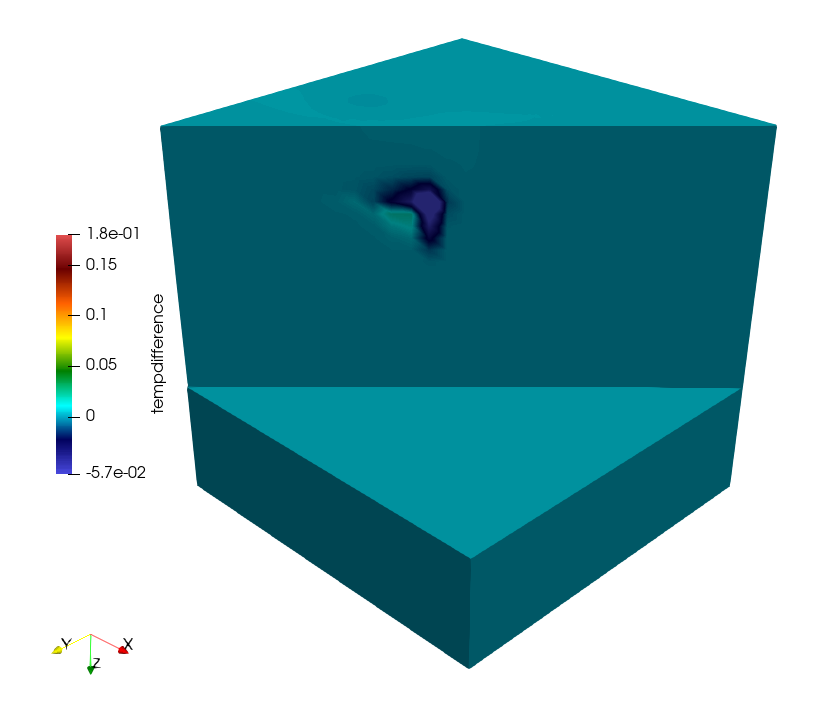}
         \caption{Difference at $t=0.25$}
         \label{fig:zdiff025}
    \end{subfigure}
    \begin{subfigure}{.5\linewidth}
         \centering
         \includegraphics[width=\textwidth]{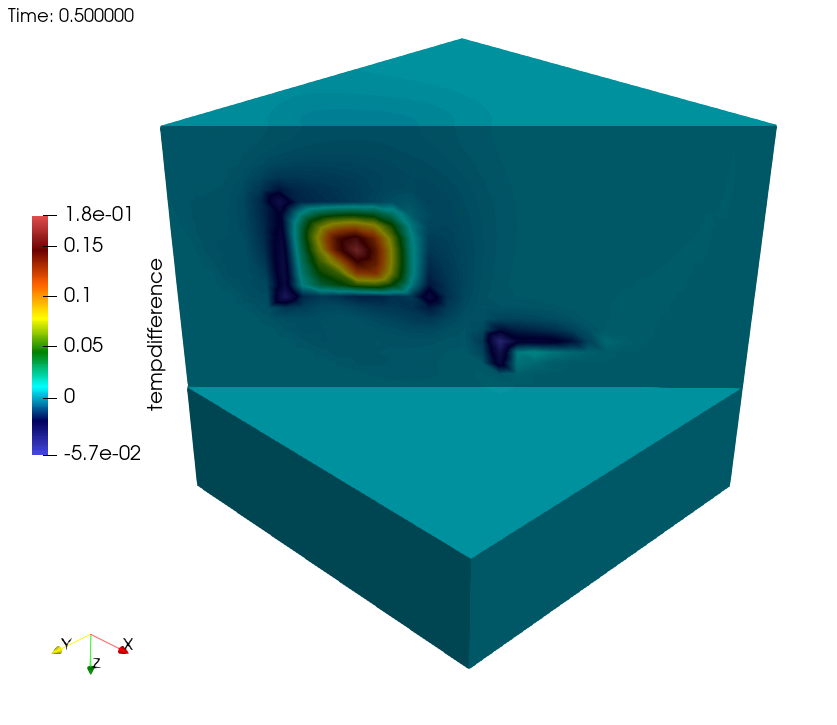}
         \caption{Difference at $t=0.5$}
         \label{fig:zdiff05}
    \end{subfigure}
    \begin{subfigure}{.5\linewidth}
         \centering
         \includegraphics[width=\textwidth]{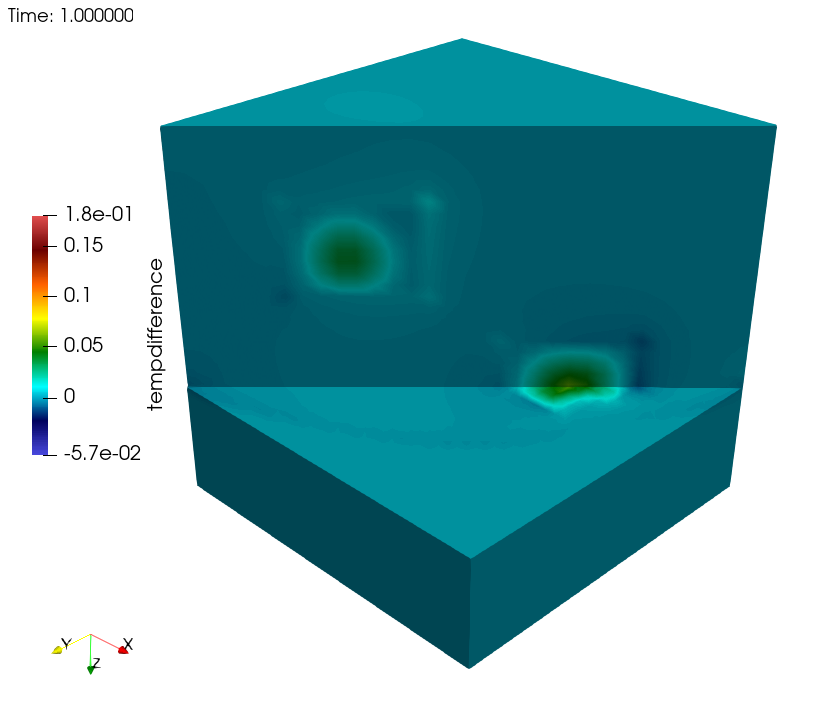}
         \caption{Difference at $t=1.0$}
         \label{fig:zdiff10}
    \end{subfigure}
    \caption{Sample Solution of the Radiative Diffusion Problem and Plots of the Difference Between the Solutions at $Z=5$ and $Z=15$ at Different Times}
    \label{fig:moosesol}
\end{figure}

In this Section, we present an application of parametric DMD to a multiphysics radiative transfer problem. Radiation transport is approximated using diffusion theory and the governing laws are given in Eqs.~\eqref{eq:radeq} (see also \cite{mousseau2003}):
\begin{subequations}\label{eq:radeq}
\begin{equation}\label{eq:radeq2}
    \frac{\partial T}{\partial t} - \nabla \cdot (D_T\nabla T) = -\sigma_a(T^4-E)\,,
\end{equation}
\begin{equation}\label{eq:radeq1}
    \frac{\partial E}{\partial t} - \nabla \cdot (D_r\nabla E) = \sigma_a(T^4-E)\,,
\end{equation}
\end{subequations}
where $E$ is the radiation energy and $T$ is the material temperature. The absorption opacity $\sigma_a$ is given by Eq.~\eqref{eq:abs} where $\alpha$ is a parameter, which we consider to be uncertain and thus can vary, and $Z$ is the atomic number of the material, which will also be a parameter that varies as the material composition is modified in the problem:
\begin{equation}\label{eq:abs}
    \sigma_a(T) = \frac{Z^\alpha}{T^\alpha}
\end{equation}
$D_t$ is the material conduction coefficient, defined in Eq.~\eqref{eq:dt}, and $D_r$ is the radiation diffusion coefficient, defined in Eq.~\eqref{eq:dr}.
\begin{subequations}
\begin{equation}\label{eq:dt}
    D_T(T) = 10^{-2} \times T^{5/2}
\end{equation}
\begin{equation}\label{eq:dr}
    D_r(E,T) = \frac{1}{3\sigma_a(T) + |\nabla E|/E}
\end{equation}
\end{subequations}
To test the various parametric DMD methods, parameters $\alpha$ and $Z$ will be allowed to vary. A 3-D version of the 2-D problem given in \cite{mousseau2003} was created. Here, the 3-D background medium has a $Z$-value of 1, while two inclusions of high-$Z$ material are present in the form of two cubes. The cubes are located at $(\tfrac{3}{32}\leq x \leq \tfrac{7}{32},\  \tfrac{9}{32}\leq y \leq \tfrac{13}{32},\ \tfrac{3}{32}\leq z \leq \tfrac{7}{32})$ and $(\tfrac{9}{32}\leq x \leq \tfrac{13}{32},\ \tfrac{3}{32}\leq y \leq \tfrac{7}{32},\ \tfrac{9}{32} \leq z \leq \tfrac{13}{32})$. Boundary conditions are reflective on the 3 adjacent faces of the domain and vacuum boundary conditions are applied on the 3 opposite faces. The initial condition for $T$ and $E$ are given in Eq.~\eqref{eq:tempic} and Eq.~\eqref{eq:radic}. 
\begin{equation}\label{eq:tempic}
    T(x,y,z,t=0) = 0.001 + 100\times\left[\exp({-100(x^2+y^2+z^2)}\right]^{\frac{1}{4}}\,,
\end{equation}
\begin{equation}\label{eq:radic}
    E(x,y,z,t=0) = 0.001+100 \times \exp({-100(x^2+y^2+z^2)}\,.
\end{equation}
The input parameter space ($\alpha$ and $Z$) is two-dimensional in this problem. A sample solution with $\alpha=3$ in the entire domain and $Z=10$ in the two inclusions is shown in \fig{fig:z10t025} at time $t=0.25$. In Figures~\ref{fig:zdiff025}-\ref{fig:zdiff10}, the difference between two solutions, with inclusion values of $Z=5$ and $Z=15$ respectively, is plotted at different time instants.

\subsubsection{Method Comparison Results}
\begin{figure}
    \begin{subfigure}{.5\linewidth}
         \centering
         \includegraphics[width=\textwidth]{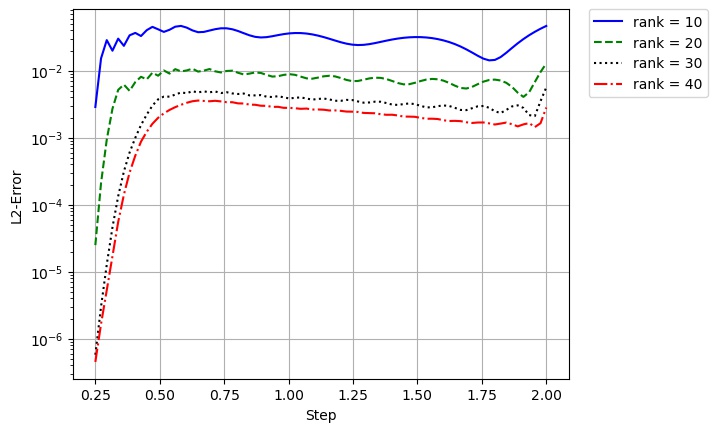}
         \caption{Stacked DMD}
         \label{fig:stackedMoose}
    \end{subfigure}%
    \begin{subfigure}{.5\linewidth}
         \centering
         \includegraphics[width=\textwidth]{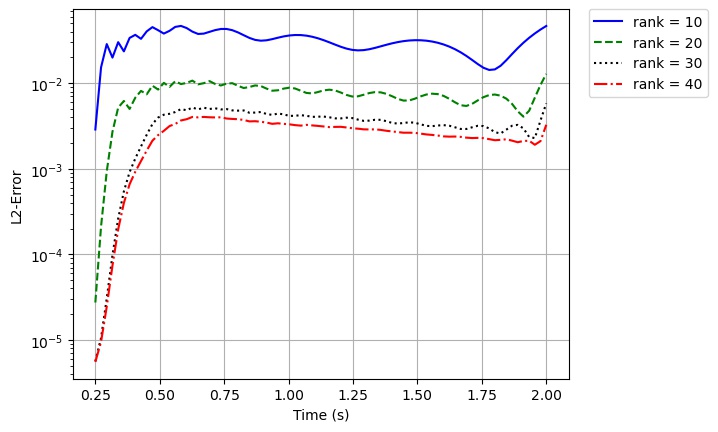}
         \caption{rKOI}
         \label{fig:koMoose}
    \end{subfigure}
    \caption{$\langle E(t_i) \rangle$ Error in Rank for two Parametric DMD Methods on the Radiative Diffusion Problem}
    \label{fig:compareMoose}
\end{figure}

In this problem, we first choose to vary only $Z$ in the two inclusion, while keeping $\alpha = 3.0$ fixed. We study varying multiple parameters in the following sections. Due to the long run times of these simulations, a set of snapshots was divided into training and testing sets. The data set includes the range $Z\in[1,15]$ at each integer value in that range. The odd-$Z$ values are used for the training snapshots and the even-$Z$ values for the testing values. In this problem we take the individual snapshots to include both the radiative energy $E$ and the temperature $T$ solution fields. The $\langle E(t_i) \rangle$ error is calculated as an average of $E$-values and $T$-values. In Fig. \ref{fig:compareMoose} we show the $\langle E(t_i) \rangle$ error of stacked DMD and rKOI as a function of rank. We also display the run time of the two methods in Table \ref{tab:Runtime}. This timing includes the SVD of the two nearest neighbors for the rKOI case and the single SVD containing all snapshots for the stacked DMD case.

\begin{table}[]
\centering
\begin{tabular}{|l|l|l|l|l|l|l|}
\hline
Rank                 & 10   & 20   & 30   & 40\\
\hline
Stacked Run Time (s) & 8.6 & 8.7 & 8.8 & 8.9 \\
\hline
rKOI Run Time (s)    & 3.4  & 3.5 & 3.5 & 3.5  \\
\hline
\% Speed Up         & 60.5\%  & 59.8\%  & 60.2\%  & 60.7\% \\
\hline
\end{tabular}
\caption{Runtime of rKOI and Stacked DMD}
\label{tab:Runtime}
\end{table}

\begin{figure}
    \centering
    \includegraphics[width=\textwidth]{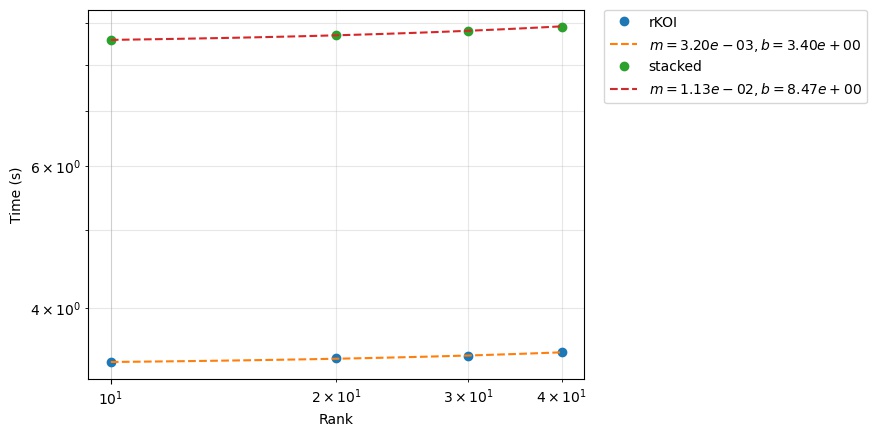}
    \caption{Rank Versus Run Time for rKOI and Stacked DMD}
    \label{fig:rankvtime}
\end{figure}

\begin{table}[]
\centering
\begin{tabular}{|l|l|l|l|l|}
\hline
Rank                 & 10   & 20   & 30   & 40 \\
\hline
rKOI Error & 3.11 E-2 & 7.49 E-3 & 3.49 E-3 & 2.65 E-3\\
\hline
Stacked Error & 3.11 E-2  & 7.50 E-3 & 3.24 E-3 & 2.16 E-3\\
\hline
\end{tabular}
\caption{$\langle \bar E  \rangle$ Error in Rank for two Parametric DMD Methods on the Radiative Diffusion Problem}
\label{tab:mooseerr}
\end{table}

In this problem we see rKOI and stacked DMD perform similarly in terms of $\langle \bar E  \rangle$ error with stacked DMD slightly beating out rKOI as shown in Table \ref{tab:mooseerr}. However rKOI sees a speed up of up to $60.5\%$ over stacked DMD, this is due to rKOI only requiring the two nearest neighbors to evaluate a point in-between while stacked DMD requires the SVD of all of the snapshots. In \ref{fig:rankvtime} we show the rank plotted versus time. The slope of this line will show how the two methods scale as a function of rank. It is clear that rKOI will scale more efficiently than stacked DMD from the respective $1.1E-3$ and $3.2E-2$ slopes. Our new method therefore is able to perform very close to the state of the art in terms of error while running up to $60\%$ faster while also scaling favorably in terms of rank. These characteristics make our method much more attractive than stacked DMD if timing is of the upmost concern. This in addition to the superior performance on problems containing advection like the Incident-Jet Problem makes our new rKOI method more attractive in a large number of scenarios.
\subsection{Parametric DMD Performance Study}

In this section we take the Radiative Diffusion problem and study the behavior of the rKOI and stacked DMD in methods more detail. We will look at the error as a function of number of snapshots, performance on a problem with 2 parameters, and how the multi-physics nature of the problem effects the behavior.

\begin{figure}
    \begin{subfigure}{.5\linewidth}
        \centering
        \includegraphics[width=\textwidth]{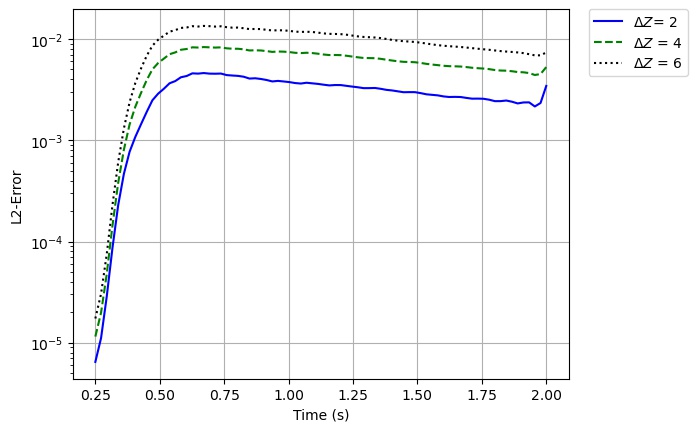}
        \caption{$\langle E(t_i) \rangle$ Error in Rank for rKOI with Varied Parameter Spacing}
        \label{fig:moosecoarserKOI}
    \end{subfigure}%
    \begin{subfigure}{.5\linewidth}
        \centering
        \includegraphics[width=\textwidth]{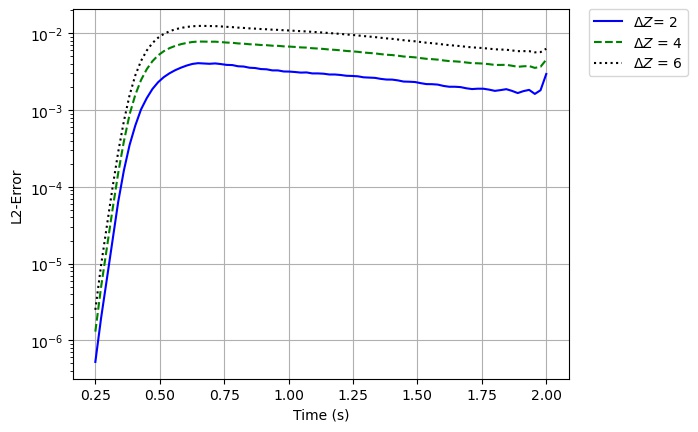}
        \caption{$\langle E(t_i) \rangle$ Error in Rank for Stacked DMD with Varied Parameter Spacing}
        \label{fig:moosecoarsestacked}
    \end{subfigure}
    \caption{$\langle E(t_i) \rangle$ Error in Rank for two Parametric DMD Methods for Varied Parameter Spacing}
    \label{fig:comparecoarse}
\end{figure}

First we will look at the $\langle E(t_i) \rangle$ error as a function of the number of snapshots. Using the same parameter range for $Z$ we use different $\Delta Z$ to generate 3 training sets, one with $\Delta Z=2$ ($Z = [1,3,5,7,9,11,13,15]$), one with $\Delta Z=4$ ($Z = [1,5,9,13]$), and one with $\Delta Z=6$ ($Z = [1,7,13]$). The testing set for this problem is the even integers from 2 to 12. We used a rank 40 reconstruction for all values of $\Delta Z$. The results are shown in Fig.~\ref{fig:comparecoarse}. 

Both methods perform similarly relative to their performance in rank. As could be expected the more snapshots that are used the lower the $\langle E(t_i) \rangle$ error of the solution is. This property could potentially be used in the future to determine how many full order runs are needed in order to achieve a certain $\langle E(t_i) \rangle$ error in the the reduced solution over the entire parameter domain. This property is also very desirable because it means that before performing DMD the more snapshots that are generated the lower the error will be ensured to be at a new parameter value.

\begin{figure}
    \begin{subfigure}{.5\linewidth}
        \centering
        \includegraphics[width=\textwidth]{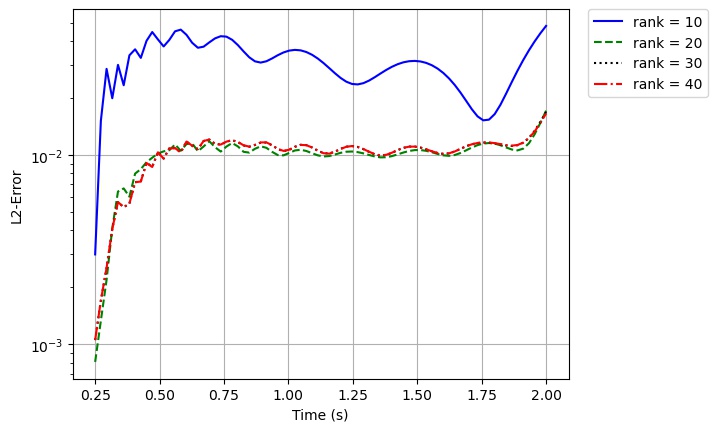}
        \caption{$\langle E(t_i) \rangle$ Error in Rank for rKOI on 2 Parameter Problem}
        \label{fig:moose2drkoi}
    \end{subfigure}%
    \begin{subfigure}{.5\linewidth}
        \centering
        \includegraphics[width=\textwidth]{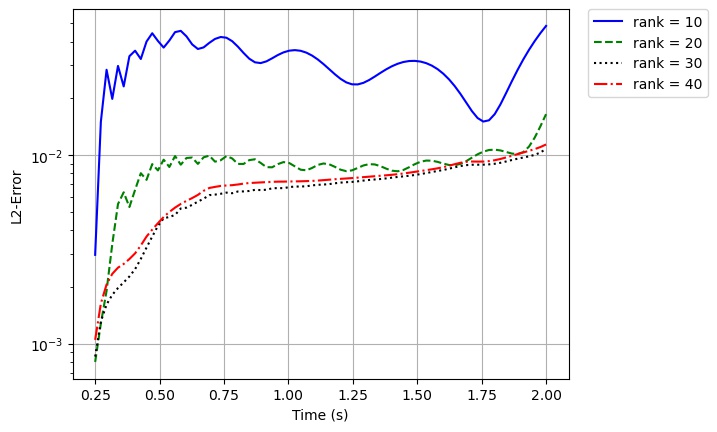}
        \caption{$\langle E(t_i) \rangle$ Error in Rank for Stacked DMD on 2 Parameter Problem}
        \label{fig:moose2dstacked}
    \end{subfigure}
    \caption{$\langle E(t_i) \rangle$ Error in Rank for two Parametric DMD Methods for 2 Parameter Problem}
    \label{fig:compare2d}
\end{figure}

Next we will expand the problem to a 2-parameter problem where in addition to varying $Z$ in the snapshots we also vary the $\alpha$ of the cubes. We will vary Z between 5 and 15 including every whole value of Z and 11 values of $\alpha$ equally spaced in $\alpha \in [2.5,3.5]$. We took all even $Z$ to be the testing set as well as every $\alpha$ with an even tenths place. The set containing each pairing of $Z$ and $\alpha$ testing points was used as the testing set and the set containing each paring of $Z$ and $\alpha$ training points was used as the training set. We evaluated at each point using the same ending criteria as before, but with only using 5 testing points instead of the entire set and the results as a function of rank are shown in Fig.~\ref{fig:compare2d}. 

Even though this problem has many more snapshots than the single parameter case the error is increased with a much higher lower error bound. This error bound differs between the two methods with the stacked DMD method having a slightly lower bound. The distances in 2-D space are larger than in 1-D parameter space making the distance between parameters for the interpolation required to get the same performance much smaller requiring more snapshots overall. This method is still able to get a reasonable error much faster than the full order model however, meaning that there is still reason to consider DMD even in multi-parameter cases.

\begin{figure}
    \begin{subfigure}{0.5\linewidth}
        \centering
        \includegraphics[width=\textwidth]{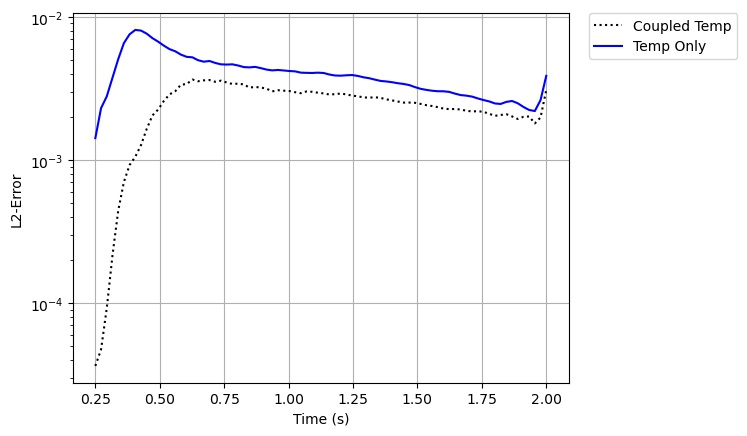}
        \caption{rKOI Coupling Effects on Temperature Solution}
        \label{fig:couplingtemprkoi}
    \end{subfigure}
        \begin{subfigure}{0.5\linewidth}
        \centering
        \includegraphics[width=\textwidth]{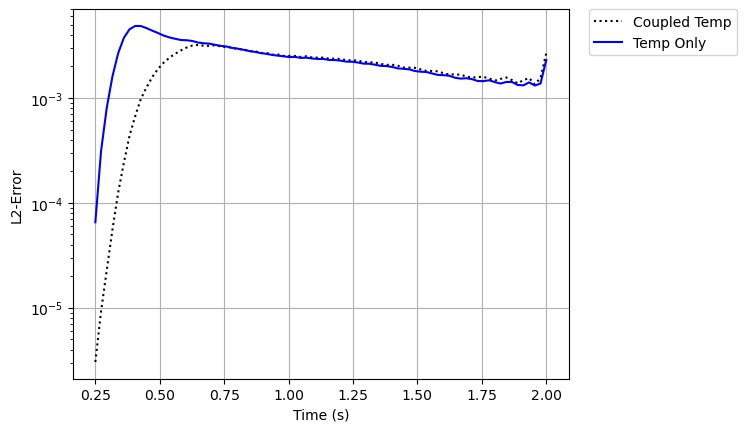}
        \caption{Stacked DMD Coupling Effects on Temperature Solution}
        \label{fig:couplingtempstacked}
    \end{subfigure}
    \begin{subfigure}{0.5\linewidth}
         \centering
         \includegraphics[width=\textwidth]{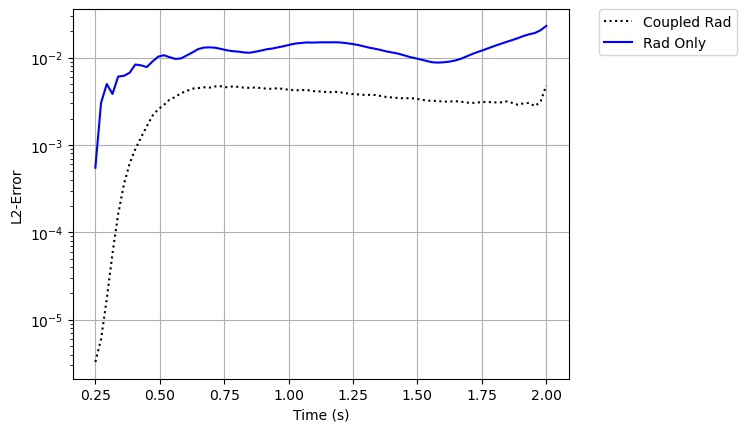}
         \caption{rKOI Coupling Effects on Radiation Solution}
         \label{fig:couplingradrkoi}
    \end{subfigure}
    \begin{subfigure}{0.5\linewidth}
         \centering
         \includegraphics[width=\textwidth]{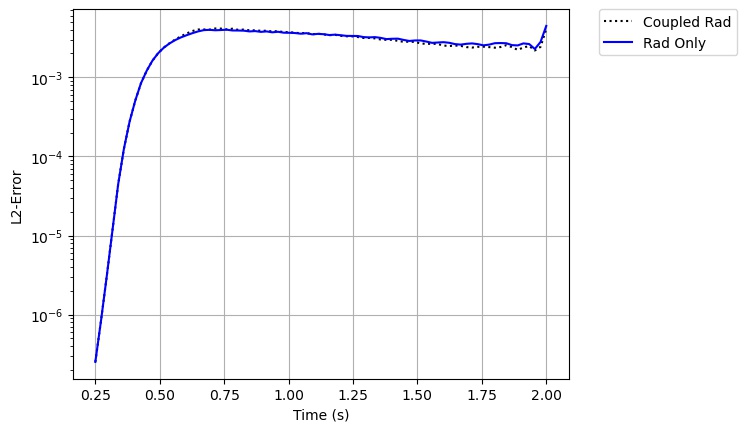}
         \caption{Stacked DMD Coupling Effects on Radiation Solution}
         \label{fig:couplingradstacked}
    \end{subfigure}
    \caption{Study of Coupling Effects on the $\langle E(t_i) \rangle$ Error of two Parametric DMD Methods}
    \label{fig:coupling}
\end{figure}

Next we will study how the coupling of this problem effects parametric DMD. As a reminder this problem is a set of coupled differential equations in both E and T. In all previous parts the error was reported averaged over both E and T. In this part we will present the error over only one portion of the solution at a time. In addition we will also perform DMD taking the snapshot to include both E and T as we have in the previous parts and perform DMD independently using only E or only T as the snapshot to predict only E or only T respectively. In this problem we go back to only vary $Z$ as the parameter and keeping $\alpha = 3.0$. We take $Z\in[1,15]$ with odd values composing the training set and even values composing the testing set. We use the same ending criteria as before and plot the results in Fig. \ref{fig:coupling}.

In the rKOI case for both E and T the coupled DMD performs better at all time steps than the DMD performed with only E or T as the snapshot. This effect is less pronounced in the stacked DMD case which means it might be preferable to use the stacked DMD method for a multi-physics problem where only one of the physics is known. Better performance in the full data case makes sense given that the coupled nature of this physics would be better captured with both halves of the snapshot instead of only one.

We have shown that both rKOI and stacked DMD have a favorable relationship between error and parameter spacing, that in a multi-parameter case both parametric DMD methods have significantly worse performance than the single parameter case with the stacked DMD slightly out performing rKOI, and that both parametric DMD methods do significantly better when they are done with the full multi-physics snapshot instead of just one physics. All of these results show that the rKOI method is similarly robust to the state of the art method.

\section{Conclusions} \label{sec:ccl}
In this article, we analyzed three methods for extending dynamic mode decomposition (DMD) to deal with parametric problems.
We analyzed the stacked-snapshots approach proposed by Sayadi et al.~\cite{sayadi2015} and we proposed two new methods based on the interpolation of the Koopman eigenpair (rEPI) and the components of the reduced Koopman operartor matrix (rKOI).
Whereas the approach proposed by Sayadi et al.~\cite{sayadi2015} scale, at best, linearly with the number of parameters in the training set, the two novel methods proposed scale only with the number of training parameters used for the interpolation set. 
This makes the methods proposed more apt for multi-parametric problems where a large number of parameter samples are needed over the training set.
Furthermore, a key limitation identified in the method proposed by Sayadi et al.~\cite{sayadi2015} is that it forces the frequency of evolution in the DMD modes to be parameter-independent.
The new methods proposed overcome this limitation, but require smoothness in the reduced Koopman operator eigenpair (rEPI) or components (rKOI) over the parameters in the interpolation set.
Note that non-linear manifold interpolation, e.g., Grassmann manifold, matrix interpolation, e.g., solving the orthogonal Procrustes problem, or polynomial regression, e.g., radial-basis function regression, may mitigate these smoothness requirements.
However, we have chosen to use simple linear interpolation methods to better showcase the performance of the parametric-DMD methods, rather than the one of the interpolation methods used.

For the three problems analyzed, we observed that the smoothness requirements in the reduced Koopman operator eigenpair is often too restrictive.
The performance of the rEPI method deteriorated at high ranks in the DMD reconstruction.
On the contrary, the smoothness assumption in the components of the reduced Koopman operator was often less constraining.
In fact, the rKOI method surpassed in accuracy the stacked DMD method in problems where the solution field evolution was strongly dependent on the uncertain parameter, while it was significantly superior in computational performance in every problem.
For all cases, the rKOI method was able to achieve a relative $L_2$ errors below 1\% along the whole problem time-span.
Moreover, for the multi-physics radiation diffusion problem, the rKOI method was also able to achieve average $L_2$ errors below 1\% when reconstructing one field at the time, independently of the other.
This single-physics reconstruction possibility offered by the rKOI method further reduces its computational burden.
Note, however, that the performance of the rKOI method in single-physics reconstruction was inferior than the one when reconstructing the whole multi-physics problem. 
This deterioration in performance was not observed when using the stacked DMD approach.

Future work should look at improving the rKOI method and rEPI methods with adapted interpolation methods, e.g., non-linear Grassmann manifold or matrix interpolation methods. Also, these framework should be applied to different problems to continue testing and improving their performance.

\section*{Acknowledgments}
This work was performed under the auspices of the U.S. Department of Energy by Lawrence Livermore National Laboratory under Contract DE-
AC52-07NA27344. The authors thank Youngsoo Choi for many fruitful discussions regarding ROM and DMD.

\bibliographystyle{unsrt}
\bibliography{main_arxiv}

\end{document}